\documentclass[11pt]{article}

\usepackage{amssymb, amsmath, amsthm, latexsym}
\usepackage{fullpage, color}
\usepackage{url}
\usepackage{algorithm}
\usepackage{algpseudocode}
\usepackage{tabularx}
\usepackage{paralist}
\usepackage{mathtools}
\usepackage{tcolorbox}
\usepackage{xcolor}

\usepackage{bbm} 

\usepackage{makecell}
\usepackage{multirow}
\usepackage{booktabs}

\usepackage[flushleft]{threeparttable} 

\usepackage{caption}
\usepackage{multirow}
\usepackage{colortbl}
\definecolor{bgcolor}{rgb}{0.8,1,1}
\definecolor{bgcolor2}{rgb}{0.8,1,0.8}
\definecolor{bgcolor3}{rgb}{0.67, 0.94, 0.82}

\usepackage{nicefrac}       

\definecolor{niceblue}{rgb}{0.0,0.19,0.56}

\usepackage{hyperref}
\hypersetup{colorlinks,linkcolor={blue},citecolor={niceblue},urlcolor={blue}}

\usepackage[numbers]{natbib}
\usepackage{sidecap}

\usepackage{pifont}
\usepackage{subfigure}

\newcommand{\R}{\mathbb{R}}

\def\<#1,#2>{\left\langle #1,#2\right\rangle}

\usepackage{mdframed} 
\usepackage{thmtools}

\definecolor{shadecolor}{gray}{0.9}
\declaretheoremstyle[
headfont=\normalfont\bfseries,
notefont=\mdseries, notebraces={(}{)},
bodyfont=\normalfont,
postheadspace=0.5em,
spaceabove=1pt,
mdframed={
  skipabove=8pt,
  skipbelow=8pt,
  hidealllines=true,
  backgroundcolor={shadecolor},
  innerleftmargin=4pt,
  innerrightmargin=4pt}
]{shaded}

\declaretheorem[style=shaded,within=section]{definition}
\declaretheorem[style=shaded,sibling=definition]{theorem}

\declaretheorem[style=shaded,sibling=definition]{corollary}

\declaretheorem[style=shaded,sibling=definition]{lemma}

\usepackage[colorinlistoftodos,bordercolor=orange,backgroundcolor=orange!20,linecolor=orange,textsize=scriptsize]{todonotes}

\newcommand{\cO}{{\cal O}}

\def\R{\mathbb{R}}

\def\R{\mathbb R}

\def\la{\langle}
\def\ra{\rangle}

\usepackage{xspace}

\newcommand{\algname}[1]{{\sf #1}\xspace}

\usepackage{hyperref}

\usepackage{accents}
\newlength{\dhatheight}

\begin{document}
\author{Marina Danilova$^{1,2}$}
\title{\textbf{On the Convergence Analysis of}\\ \textbf{Aggregated Heavy-Ball Method}\thanks{The research was supported by Russian Foundation
for Basic Research (Theorem~\ref{thm:AggHB_non_convex}, project No.\ 20-31-90073) and by Russian Science Foundation (Theorem~\ref{thm:AggHB_main_result}, project No.\
21-71-30005).}}
\date{$^1$Institute of Control Sciences of Russian Academy of Sciences, Russia\\ $^2$Moscow Institute of Physics and Technology, Russia}
\maketitle
\begin{abstract}
Momentum first-order optimization methods are the workhorses in various optimization tasks, e.g., in the training of deep neural networks. Recently, Lucas et al.\ (2019) \cite{lucas2019aggregated} proposed a method called Aggregated Heavy-Ball (\algname{AggHB}) that uses multiple momentum vectors corresponding to different momentum parameters and averages these vectors to compute the update direction at each iteration. Lucas et al.\ (2019) \cite{lucas2019aggregated} show that \algname{AggHB} is more stable than the classical Heavy-Ball method even with large momentum parameters and performs well in practice. However, the method was analyzed only for quadratic objectives and for online optimization tasks under uniformly bounded gradients assumption, which is not satisfied for many practically important problems. In this work, we address this issue and propose the first analysis of \algname{AggHB} for smooth objective functions in non-convex, convex, and strongly convex cases without additional restrictive assumptions. Our complexity results match the best-known ones for the Heavy-Ball method. We also illustrate the efficiency of \algname{AggHB} numerically on several non-convex and convex problems.
\end{abstract}


\section{Introduction}\label{sec:intro}

Momentum \cite{polyak1964some} and acceleration \cite{nesterov1983method} are popular techniques for speeding up first-order optimization methods both from practical and theoretical perspectives. Historically, one of the first examples of such methods is Heavy-Ball (\algname{HB}) method proposed by B.\ Polyak in 1964 \cite{polyak1964some}. This method received a lot of attention from various research communities due to its efficiency in different convex and, more importantly, non-convex problems \cite{danilova2020recent}. In particular, during the last few years a lot of variants of \algname{HB} were proposed and analyzed by machine learning (ML) researchers, especially due to its efficiency in computer vision tasks \cite{sutskever2013importance}.

Recently, another modification of \algname{HB} called Aggregated Heavy-Ball (\algname{AggHB}) method was proposed in \cite{lucas2019aggregated}. In contrast to \algname{HB}, \algname{AggHB} has $m \geq 1$ different momentum parameters and $m$ corresponding momentum vectors. An average of these vectors is used as an update direction at each iteration. Such an averaging helps to make the method more stable via reducing the oscillations of the iterates, as the authors of \cite{lucas2019aggregated} illustrated empirically. Moreover, the numerical results from \cite{lucas2019aggregated} show the superiority of \algname{AggHB} to \algname{HB} at training several ML models.

\subsection{Motivational Example}
\begin{figure}[t]
    \centering
    \includegraphics[width=0.41\textwidth]{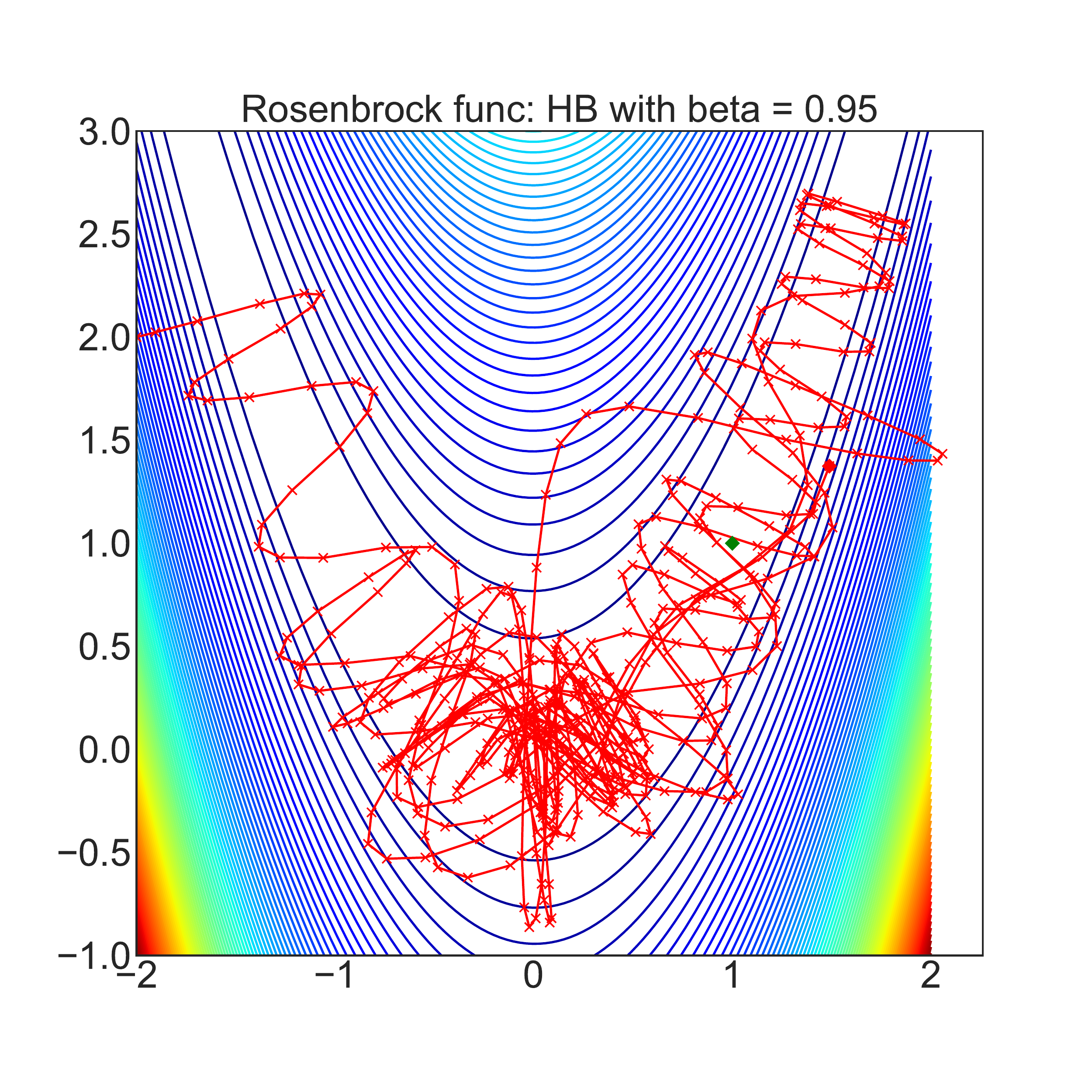}
    \includegraphics[width=0.41\textwidth]{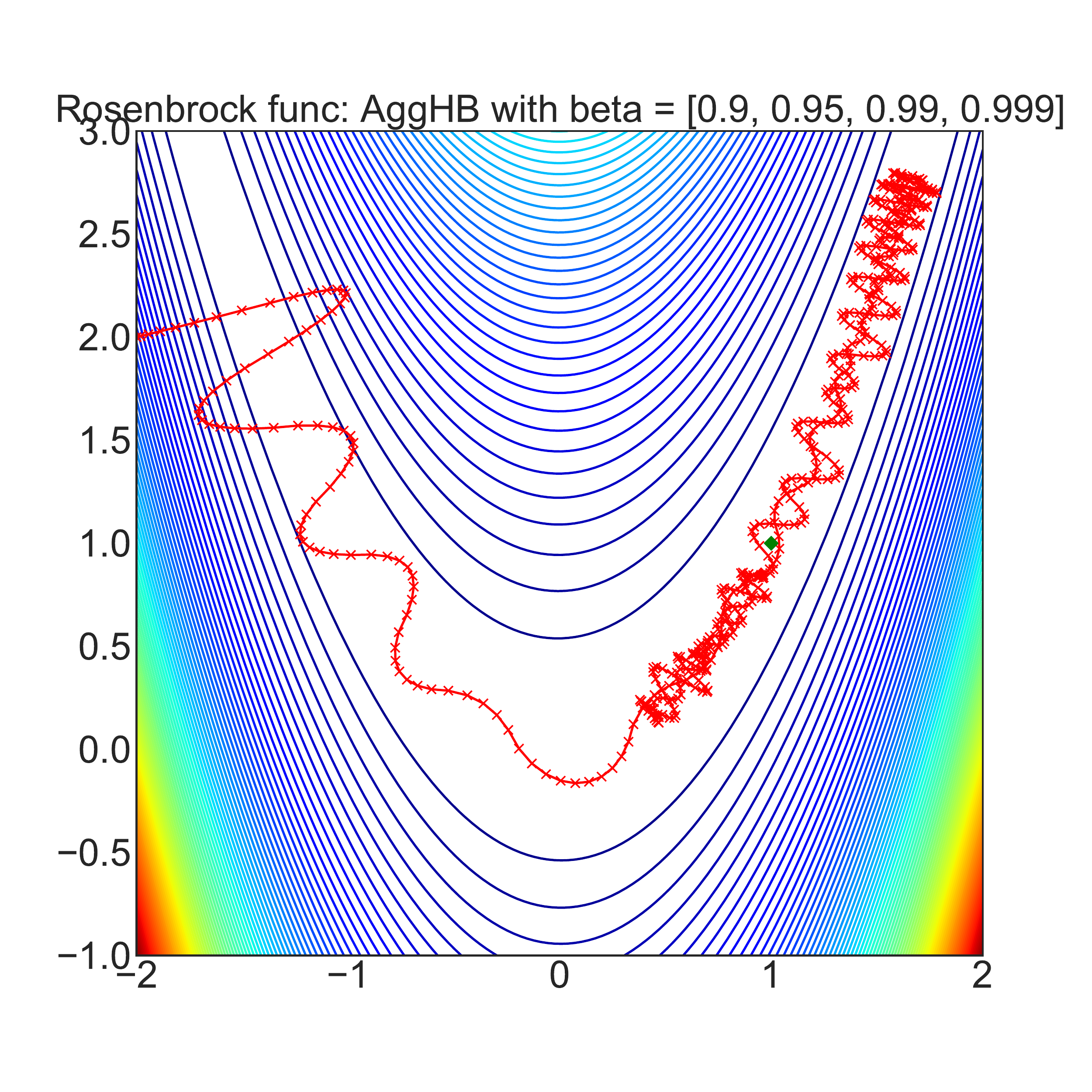}
    \caption{Trajectories of \algname{HB} (left) and \algname{AggHB} (right) with different momentum parameters $\beta$ applied to minimize Rosenbrock function. Stepsize $\gamma$ was tuned for each method. We use the package from \cite{Novik_torchoptimizers} for the visualization.}
    \label{fig:rosenbrock}
\end{figure}

In this section, we consider the behavior of \algname{AggHB} on Rosenbrock function, which is well-known non-convex test functions. The set of momentum parameters for \algname{AggHB} were chosen as $[0.9, 0.95, 0.99, 0.999]$ (see Algorithm~\ref{alg:AggHB}) and for \algname{HB} a standard momentum parameter $\beta = 0.95$ was taken (see Algorithm~\ref{alg:HB_m}). Stepsize $\gamma$ was tuned for each method. The results are presented in Figure~\ref{fig:rosenbrock}. We observe much smaller oscillations for \algname{AggHB} than for \algname{HB}. Moreover, the trajectory of  \algname{AggHB} achieves better accuracy. This example motivates the detailed study of \algname{AggHB} and, in particular, the theoretical study of its convergence.

\subsection{Our Contributions}

However, a little is known about theoretical convergence guarantees for \algname{AggHB}. In particular, the authors of \cite{lucas2019aggregated} analyzed \algname{AggHB} for quadratic optimization problems, which is a very small class of problems, and for convex online optimization problems such that the gradients of the objective function are bounded on the whole domain. The former assumption is not satisfied for many practically important tasks. \emph{In this paper, we remove this limitation and derive new convergence results for \algname{AggHB} for smooth non-convex and (strongly) convex problems.}

Our main contributions can be summarized as follows.

\begin{itemize}
    \item[$\diamond$] \textbf{First analysis of \algname{AggHB} for non-convex problems.} For the problems with smooth but not necessary convex objective function $f$, we prove that \algname{AggHB} finds an $\varepsilon$-stationary point (point $x$ such that $\|\nabla f(x)\| \leq \varepsilon$) after $\cO(\nicefrac{1}{\varepsilon^2})$ iterations neglecting the dependence on momentum parameters, smoothness constant, and initial functional suboptimality. When $m=1$ we recover the complexity of \algname{HB} and when $m > 1$ our rate is better than the corresponding rate of \algname{HB} with maximal momentum parameter (see Theorem~\ref{thm:AggHB_non_convex} and Corollary~\ref{cor:AggHB_complexity_non_cvx} for the details).
    \item[$\diamond$] \textbf{First analysis of \algname{AggHB} without bounded gradient assumption.} In the smooth (strongly) convex case, we derive the first complexity upper bounds for \algname{AggHB} without assuming that the gradients are uniformly bounded. As in the non-convex case, we recover the complexity of \algname{HB} when $m=1$ and our rate is better than the corresponding rate of \algname{HB} with maximal momentum parameter when $m > 1$ (see Theorem~\ref{thm:AggHB_main_result} and Corollary~\ref{cor:AggHB_complexity}).
    \item[$\diamond$] \textbf{Numerical experiments.} We compare the performance of \algname{AggHB} and \algname{HB} on the logistic regression problem with $\ell_2$-regularization and special non-convex regularization. In our experiments, \algname{AggHB} converges faster than \algname{HB}.
\end{itemize}

\subsection{Technical Preliminaries}
We consider an unconstrained minimization problem
\begin{equation}
    \min\limits_{x\in \R^n}f(x),\label{eq:main_problem}
\end{equation}
where function $f: \R^n \to \R$ is $L$-smooth, i.e., for all $x,y\in \R^n$
\begin{equation}
    \|\nabla f(x) - \nabla f(y)\|_2 \le L\|x - y\|_2. \label{eq:L_smooth_def}
\end{equation}
Next, we assume that $f(x)$ is either bounded from below $f_{\inf} = \inf_{x\in\R^n} f(x) > -\infty$ or $\mu$-strongly convex
\begin{equation}
    f(y) \ge f(x) + \langle \nabla f(x), y-x \rangle + \frac{\mu}{2}\|y-x\|_2^2. \label{eq:mu_str_cvx_def}
\end{equation}
The notation we use is standard for optimization literature \cite{polyak1987introduction,nesterov2018lectures}, e.g., by $x_*$ we denote the solution of \eqref{eq:main_problem}, the distance from the starting point to the solution is denoted by $R_0 = \|x_0 - x_*\|_2$.

\subsection{Related Work}

\begin{algorithm}[t]
\caption{Heavy-Ball method (\algname{HB})}
\label{alg:HB_m}   
\begin{algorithmic}[1]
\Require starting points $x_0$, $x_1$ (by default $x_0 = x_1$), number of iterations $N$, stepsize $\gamma > 0$, momentum parameter $\beta \in [0,1]$
\For{$k=0,\ldots, N-1$}
\State  $V_k=\beta V_{k-1} + \nabla f(x_k)$
\State $x_{k+1} = x_k - \gamma V_k$
\EndFor
\Ensure $x_N$ 
\end{algorithmic}
\end{algorithm}

\paragraph{Theoretical convergence guarantees for \algname{HB}.} The first convergence analysis of Heavy-Ball method (\algname{HB}, Algorithm~\ref{alg:HB_m}) was given in the original work by B.\ Polyak in 1964 \cite{polyak1964some}, where \emph{local}  $\cO(\sqrt{\nicefrac{L}{\mu}}\log(\nicefrac{1}{\varepsilon}))$ convergence rate was shown for twice continuously differentiable $L$-smooth and $\mu$-strongly convex functions. After 50 years Ghadimi et al.\ (2015) \cite{ghadimi2015global} derived the first \emph{global} convergence rates for \algname{HB} (and its version with averaging). In particular, they shown $\cO(\nicefrac{L}{\mu}\log(\nicefrac{1}{\varepsilon}))$ and $\cO\left(\nicefrac{LR_0^2}{\varepsilon}\right)$ complexity bounds for $L$-smooth $\mu$-strongly convex and convex functions respectively. In contrast to the local convergence guarantees, these rates are not accelerated \cite{nemirovsky1983problem,nesterov1983method}. Although one can improve the analysis of \algname{HB} for quadratic functions and get asymptotically accelerated rate \cite{lessard2016analysis}, it is still unclear whether this result can be generalized to the general non-quadratic functions. The non-triviality of this question is supported by the negative result from \cite{taylor2019stochastic} showing that one cannot derive accelerated rate of \algname{HB} for the standard choice of parameters using quadratic potentials in the analysis.

\paragraph{\algname{HB} with aggregation and averaging.} As we we already mentioned, Aggregated Heavy-Ball method (\algname{AggHB}, Algorithm~\ref{alg:AggHB}) was proposed in \cite{lucas2019aggregated}, where authors empirically shown that aggregation helps to stabilize the methods behavior, speeds up the method in practice, and they also derive some convergence guarantees under uniformly bounded gradients assumption in the stochastic case. Recently, in \cite{danilova2021averaged}, another approach for stabilizing \algname{HB} was considered. In particular, the authors of \cite{danilova2021averaged} considered several averaging techniques for \algname{HB} and shown that they help to reduce the maximal deviation of the method and improve the performance of the method in practice.

\section{Analysis of Aggregated Heavy-Ball Method}\label{sec:convergence_guarantees}

In this section we propose a new convergence analysis for Aggregated Heavy-Ball method (\algname{AggHB}, Algorithm~\ref{alg:AggHB}). The key difference between \algname{HB} and \algname{AggHB} is that instead of one direction determined by parameter $\beta$ the method uses to the vector of ,o,entum parameters $\beta = [\beta_1, \dots, \beta_m]$ and takes and average over $m$ corresponding directions. When $m = 1$ \algname{AggHB} recovers \algname{HB}. Moreover, we consider a slight generalization of the method proposed in \cite{lucas2019aggregated}, since we allow to use different stepsizes for different momentum parameters.

Following \cite{mania2017perturbed,yang2016unified} we consider \emph{perturbed/virtual} iterates:
\begin{equation}
    \widetilde{x}_k = x_k - \frac{1}{m}\sum\limits_{i=1}^{m}\frac{\beta_i \gamma_i}{1 - \beta_i}V_{k-1}^{(i)},\; k\ge 0. \label{eq:virtual_iterates_AggHB}
\end{equation}
This representation is used for the analysis only and there is no need to compute this sequence when running the method. Virtual iterates satisfy the following useful recursion: for all $k \ge 0$
\begin{eqnarray}
    \widetilde{x}_{k+1} &=& x_{k+1} - \frac{1}{m}\sum\limits_{i=1}^{m}\frac{\beta_i \gamma_i}{1 - \beta_i}V_{k}^{(i)} = x_k -\frac{1}{m}\sum\limits_{i=1}^m \gamma_i V_k^{(i)} - \frac{1}{m}\sum\limits_{i=1}^{m}\frac{\beta_i \gamma_i}{1 - \beta_i}V_{k}^{(i)} \notag\\
    &=& x_k - \frac{1}{m}\sum\limits_{i=1}^m \frac{\gamma_i V_k^{(i)}}{1-\beta_i} = x_k - \frac{1}{m}\sum\limits_{i=1}^{m}\frac{\beta_i \gamma_i}{1 - \beta_i}V_{k-1}^{(i)} - \frac{1}{m}\sum\limits_{i=1}^m \frac{\gamma_i}{1-\beta_i}\nabla f(x_k)\notag\\
    &=&  \widetilde{x}_{k} - \frac{1}{m}\sum\limits_{i=1}^m \frac{\gamma_i}{1-\beta_i}\nabla f(x_k). \label{eq:virtual_iter_recurrence}
\end{eqnarray}

\begin{algorithm}[t]
\caption{Aggregated Heavy-Ball method (\algname{AggHB})}
\label{alg:AggHB}   
\begin{algorithmic}[1]
\Require number of iterations $N$, stepsize $\gamma_i > 0$, momentum parameters $\{\beta_i\}_{i=1}^m \in [0,1]$, starting points $x_0$, $x_1$ (by default $x_1 = x_0 - \alpha \nabla f(x_0)$)
\For{$k=1,\ldots, N-1$}
\State $V_k^{(i)} = \beta_iV_{k-1}^{(i)} + \nabla f(x_k)$ for $i = 1,\ldots,m$
\State $x_{k+1}=x_k - \frac{1}{m}\sum\limits_{i=1}^m \gamma_{i}V_k^{(i)}$
\EndFor
\Ensure $x_N$ 
\end{algorithmic}
\end{algorithm}

\subsection{Non-Convex Case}
Below we present our main convergence result\footnote{We defer all the proofs to the Appendix.} for non-convex problems.

\begin{theorem}\label{thm:AggHB_non_convex}
    Let be $f$ is $L$-smooth and possibly non-convex function with values lower bounded by $f_{\inf}$. Assume that  
    \begin{equation}
        -\frac{A}{2}\left( 1 - \frac{CDEL^2}{2m^2} - LA\right) < 0,
        \label{eq:condition_non_cvx}
    \end{equation}
    where
    \begin{eqnarray}
        A = \frac{1}{m}\sum\limits_{i=1}^m\frac{\beta_i \gamma_{i}}{1-\beta_i},\; C = \sum\limits_{i=1}^m \frac{\gamma_{i}}{(1-\beta_i)^2}, \; D = \max\limits_{i=1,m} \frac{\gamma_{i}}{1-\beta_i},\; E = \sum\limits_{i=1}^m \frac{1}{1-\beta_i}.
        \label{eq:constants_non_cvx}
    \end{eqnarray}
    Then, for all $K\ge 1$ we have
    \begin{equation}
        \min\limits_{k=1,K} \|\nabla f(x_{k})\|^2_2 \le \frac{2}{K} \frac{f(x_0)-f_{\inf}}{A\left(1 - \frac{CDEL^2}{m^2} - LA\right)}.\label{eq:main_result_non_cvx}
    \end{equation}
\end{theorem}

The above result provides a convergence guarantee in the general non-convex case and allows to use different $\gamma_i$ such that \eqref{eq:condition_non_cvx} holds. To illustrate this result and, in particular, condition \eqref{eq:condition_non_cvx} we derive the following corollary of Theorem~\ref{thm:AggHB_non_convex}.

\begin{corollary}\label{cor:AggHB_complexity_non_cvx}
    Let the assumptions of Theorem~\ref{thm:AggHB_non_convex} hold. Assume that the stepsize is constant $\gamma_i \equiv \gamma$ for $i = 1,\ldots,m$ and consider new constants $\widetilde{\beta}$ and $\hat{\beta}$ satisfying the following conditions: $\tfrac{1}{m}\sum_{i=1}^m \tfrac{\beta_i}{(1-\beta_i)^2} = \tfrac{\widetilde{\beta}}{(1-\widetilde{\beta})^2},$ $ \tfrac{1}{m}\sum_{i=1}^m \tfrac{1}{1-\beta_i} = \tfrac{1}{1-\hat{\beta}}.$ Let
    \begin{equation*}
        \gamma = \frac{1}{L\left(\frac{2\hat{\beta}}{1-\hat{\beta}} + \sqrt{2\left(\frac{\widetilde{\beta}}{(1-\widetilde{\beta})^2} + \frac{1}{1-\hat{\beta}}\right)\frac{1}{\left(1-\max\limits_{i=1,m}\beta_i\right)(1-\hat{\beta})}}\right)}.
    \end{equation*}
    Then, to achieve $\min\limits_{k=1,K}\|\nabla f(x_k)\|^2_2 \le \varepsilon^2$ for $\varepsilon > 0$ \algname{AggHB} requires
    \begin{eqnarray}
        \cO\left(\frac{L(f(x_0) - f_{\inf})}{\varepsilon^2} + \frac{L(f(x_0) - f_{\inf}) \sqrt{\left(\frac{\widetilde{\beta}(1-\hat\beta)}{(1-\widetilde{\beta})^2} + 1\right)\frac{1}{\left(1-\max\limits_{i=1,m}\beta_i\right)\hat\beta^2}}}{\varepsilon^2}\right). \label{eq:AggHB_compl_str_non_cvx}
    \end{eqnarray}
\end{corollary}

First of all, when $m = 1$, we have $\beta = \widetilde{\beta} = \hat \beta = \max_{i=1,m}\beta_i$ and the above convergence rate can be simplified to
\begin{equation*}
    \cO\left(\frac{L(f(x_0) - f_{\inf})}{\varepsilon^2} + \frac{L(f(x_0) - f_{\inf})}{\varepsilon^2\beta(1-\beta)}\right)
\end{equation*}
that matches the rate of \algname{HB} in the non-convex case (e.g., see \cite{defazio2020momentum}). Next, constants $\widetilde{\beta}$ and $\hat\beta$ can be viewed as special ``averaged'' momentum parameters. Indeed, we know that
\begin{eqnarray*}
    \frac{\min_{i=1,m}\beta_i}{(1-\min_{i=1,m}\beta_i)^2} &\leq \frac{1}{m}\sum\limits_{i=1}^m \frac{\beta_i}{(1-\beta_i)^2} &\leq \frac{\max_{i=1,m}\beta_i}{(1-\max_{i=1,m}\beta_i)^2},\\
    \frac{1}{1-\min_{i=1,m} \beta_i} &\leq \frac{1}{m}\sum\limits_{i=1}^m \frac{1}{1-\beta_i} &\leq \frac{1}{1-\max_{i=1,m} \beta_i},
\end{eqnarray*}
since $\tfrac{x}{(1-x)^2}$ and $\tfrac{1}{1-x}$ are increasing functions for $x\in (0,1)$, i.e., $\widetilde{\beta}, \hat \beta$ lie in $[\min_{i=1,m}\beta_i, \max_{i=1,m}\beta_i]$. This allows to use larger stepsize than maximal possible stepsize for \algname{HB} with $\beta = \max_{i=1,m} \beta_i$, i.e., the rate of \algname{AggHB} is better than the one of \algname{HB} with $\beta = \max_{i=1,m} \beta_i$.

\subsection{Convex and Strongly-Convex Cases}

\begin{lemma}\label{lem:one_iter_progress_AggHB}
    Let be $f$ is $L$-smooth and $\mu$-strongly convex. Let $\gamma_i$ and $\beta_i$ satisfy $\gamma_i > 0,$ $\beta_i \in [0,1)$, and
    \begin{equation}
       F = \frac{1}{m}\sum\limits_{i=1}^{m}\frac{\gamma_i}{1-\beta_i} \le \frac{1}{4L}.  \label{eq:constants_AggHB}
    \end{equation}
    Then, for all $k\ge 0$
    \begin{equation}
        \frac{F}{2}\!\left(f(x_k)\! -\! f(x_*)\right) \le \left(1\!-\!\frac{F\mu}{2}\right)\|\widetilde{x}_k\! -\! x_*\|^2_2 - \|\widetilde{x}_{k+1}\! -\! x_*\|^2_2 + 3LF\|x_k\! -\! \widetilde{x}_k\|^2_2. \label{eq:one_iter_progress_AggHB}
    \end{equation}
\end{lemma}

Next, it is sufficient to sum up \eqref{eq:one_iter_progress_AggHB} for $k = 0,1,\dots K$ with weights $w_k = (1 - \nicefrac{\mu F}{2})^{-(k+1)},$ $W_k = \sum_{k=0}^{K}w_k$ to get the bound on $f(\overline{x}_K) - f(x_*)$, where $\overline{x}_K = \tfrac{1}{W_K}\sum_{i=1}^K w_k(f(x_k)-f(x_*))$. To get final result one needs to upper bound the sum $3LF \sum_{k=0}^K w_k \|x_k - \widetilde{x}_k\|^2_2$. For this we consider the following lemma.
\begin{lemma}\label{lem:weighted_sum_of_momentums}
    Assume that $f$ is $L$-smooth and $\mu$-strongly convex. Let $\gamma_i$ and $\beta_i$ satisfy
    \begin{gather}
        0 < \gamma_i \le \frac{\left(1-\max_{i=1,m}\beta_i\right)(1-\beta_i)}{2\mu}, \quad \beta_i \in [0,1),  \label{eq:params_AggHB_2}\\
        F = \frac{1}{m}\sum\limits_{i=1}^{m}\frac{\gamma_i}{1-\beta_i} \le \frac{1}{4L},  \quad BF \le \frac{1 - \max_{i=1,m}\beta_i}{48 L^2},  \label{eq:constants_AggHB_2}
    \end{gather}
    where $B = \tfrac{1}{m}\sum_{i=1}^m \tfrac{\beta_i \gamma_i \left(1-\beta_i^{K+1}\right)}{(1-\beta_i)^2}$. Then, for all $k\ge 0$ and $w_k = \left(1 - \nicefrac{\mu F}{2}\right)^{-(k+1)}$
    \begin{equation}
        3LF \sum\limits_{k=0}^K w_k \|x_k - \widetilde{x}_k\|^2_2 \le \frac{F}{4}\sum\limits_{k=0}^{K}w_k\left(f(x_k) - f(x_*)\right) \label{eq:weighted_sum_of_momentums}
    \end{equation}
\end{lemma}

Combining these lemmas, we get the main result in (strongly) convex case.
\begin{theorem}\label{thm:AggHB_main_result}
    Assume that $f$ is $L$-smooth and $\mu$-strongly convex. Let $\gamma_i$ and $\beta_i$ satisfy conditions from \eqref{eq:params_AggHB_2} and \eqref{eq:constants_AggHB_2}. Then, after $K \ge 0$ iterations of \algname{AggHB} we have
    \begin{equation}
        f(\overline{x}_K) - f(x_*) \le \frac{4\|x_0 - x_*\|_2^2}{F W_K},\quad \overline{x}_K = \frac{1}{W_K}\sum\limits{i=1}^K w_k(f(x_k)-f(x_*)) \label{eq:AggHB_main_result}
    \end{equation}
    where $w_k = (1 - \nicefrac{\mu F}{2})^{-(k+1)},$ $W_k = \sum_{k=0}^{K}w_k$, i.e., 
    \begin{eqnarray}
        f(\overline{x}_K) - f(x_*) &\le& \left(1 - \frac{\mu F}{2}\right)^K\frac{4\|x_0 - x_*\|_2^2}{F},\quad \text{if } \mu > 0,\label{eq:AggHB_main_result_str_cvx}\\
        f(\overline{x}_K) - f(x_*) &\le& \frac{4\|x_0 - x_*\|_2^2}{F K},\hspace{2.37cm} \text{if } \mu = 0. \label{eq:AggHB_main_result_cvx}
    \end{eqnarray}
\end{theorem}

As in the non-convex case, the above result gives convergence guarantees in the general convex and strongly convex cases and allows to use different $\gamma_i$ such that \eqref{eq:params_AggHB_2} and \eqref{eq:constants_AggHB_2} hold. To illustrate this result and, in particular, conditions \eqref{eq:params_AggHB_2} and \eqref{eq:constants_AggHB_2} we derive the following corollary of Theorem~\ref{thm:AggHB_main_result}.

\begin{corollary}\label{cor:AggHB_complexity}
    Let the assumptions of Theorem~\ref{thm:AggHB_main_result} hold. Assume that the stepsize is constant $\gamma_i \equiv \gamma$ for $i = 1,\ldots,m$ and consider constants $\widetilde{\beta}$ and $\hat{\beta}$ satisfying the following conditions: $\tfrac{1}{m}\sum_{i=1}^m \tfrac{\beta_i}{(1-\beta_i)^2} = \tfrac{\widetilde{\beta}}{(1-\widetilde{\beta})^2},$ $\tfrac{1}{m}\sum_{i=1}^m \tfrac{1}{1-\beta_i} = \tfrac{1}{1-\hat{\beta}}$. Let
    \begin{equation*}
        \gamma = \min\left\{\frac{\left(1-\max\limits_{i=1,m}\beta_i\right)^2}{2\mu}, \frac{1-\hat{\beta}}{4L}, \frac{(1-\widetilde{\beta})\sqrt{(1-\hat{\beta})\left(1-\max\limits_{i=1,m}\beta_i\right)}}{4\sqrt{3}L\sqrt{\widetilde{\beta}}} \right\}.
    \end{equation*}
    Then, to achieve $f(\overline{x}_K) - f(x_*) \le \varepsilon$ for $\varepsilon > 0$ \algname{AggHB} requires
    \begin{eqnarray}
        \cO\Bigg(\Bigg(\frac{L}{\mu} + \frac{1-\hat{\beta}}{\left(1-\max\limits_{i=1,m}\beta_i\right)^2} + \frac{L\sqrt{\widetilde{\beta}(1-\hat{\beta})}}{\mu(1-\widetilde{\beta})\sqrt{1-\max\limits_{i=1,m}\beta_i}}\Bigg)\qquad\qquad\quad \notag\\
        \qquad\qquad\quad\cdot \ln \Bigg( \frac{R_0^2}{\varepsilon}\cdot \Bigg( L + \frac{1-\hat{\beta}}{(1-\max\limits_{i=1,m} \beta_i)^2} + \frac{L \sqrt{\widetilde{\beta}(1-\hat{\beta})}}{(1-\widetilde{\beta})\sqrt{1-\max\limits_{i=1,m} \beta_i}}\Bigg)\Bigg) \Bigg) \label{eq:AggHB_compl_str_cvx}
    \end{eqnarray}
    iterations when $\mu > 0$, and
    \begin{equation}
        \cO\left(\frac{LR_0^2}{\varepsilon} + \frac{LR_0^2\sqrt{\widetilde{\beta}(1-\hat{\beta})}}{\varepsilon(1-\widetilde{\beta})\sqrt{1-\max\limits_{i=1,m}\beta_i}}\right) \label{eq:AggHB_compl_cvx}
    \end{equation}
    iterations when $\mu = 0$, where $R_0 = \|x_0 - x_*\|_2$.
\end{corollary}

First of all, when $m = 1$, we have $\beta = \widetilde{\beta} = \hat \beta = \max_{i=1,m}\beta_i$ and the above convergence rates can be simplified to
\begin{gather*}
    \cO\left(\left(\frac{L}{\mu} + \frac{L\sqrt{\beta}}{\mu(1-\beta)}\right) \log \left( \frac{R_0^2}{\varepsilon}\cdot \left( L + \frac{L \sqrt{\beta}}{1-\beta}\right)\right) \right), \quad \text{when } \mu > 0,\\
    \cO\left(\frac{LR_0^2}{\varepsilon} + \frac{LR_0^2\sqrt{\beta}}{\varepsilon(1-\beta)}\right), \quad \text{when } \mu = 0
\end{gather*}
that matches the rate of \algname{HB} in the strongly convex and convex cases (e.g., see \cite{ghadimi2015global}). Next, as we already mentioned before, constants $\widetilde{\beta}$ and $\hat\beta$ can be viewed as special ``averaged'' momentum parameters. This allows to use larger stepsize than maximal possible stepsize for \algname{HB} with $\beta = \max_{i=1,m} \beta_i$, i.e., the rate of \algname{AggHB} is better than the one of \algname{HB} with $\beta = \max_{i=1,m} \beta_i$.

\section{Numerical Experiments}\label{sec:numerical_exp}
\begin{figure}[t]
    \centering
    \includegraphics[width=0.325\textwidth]{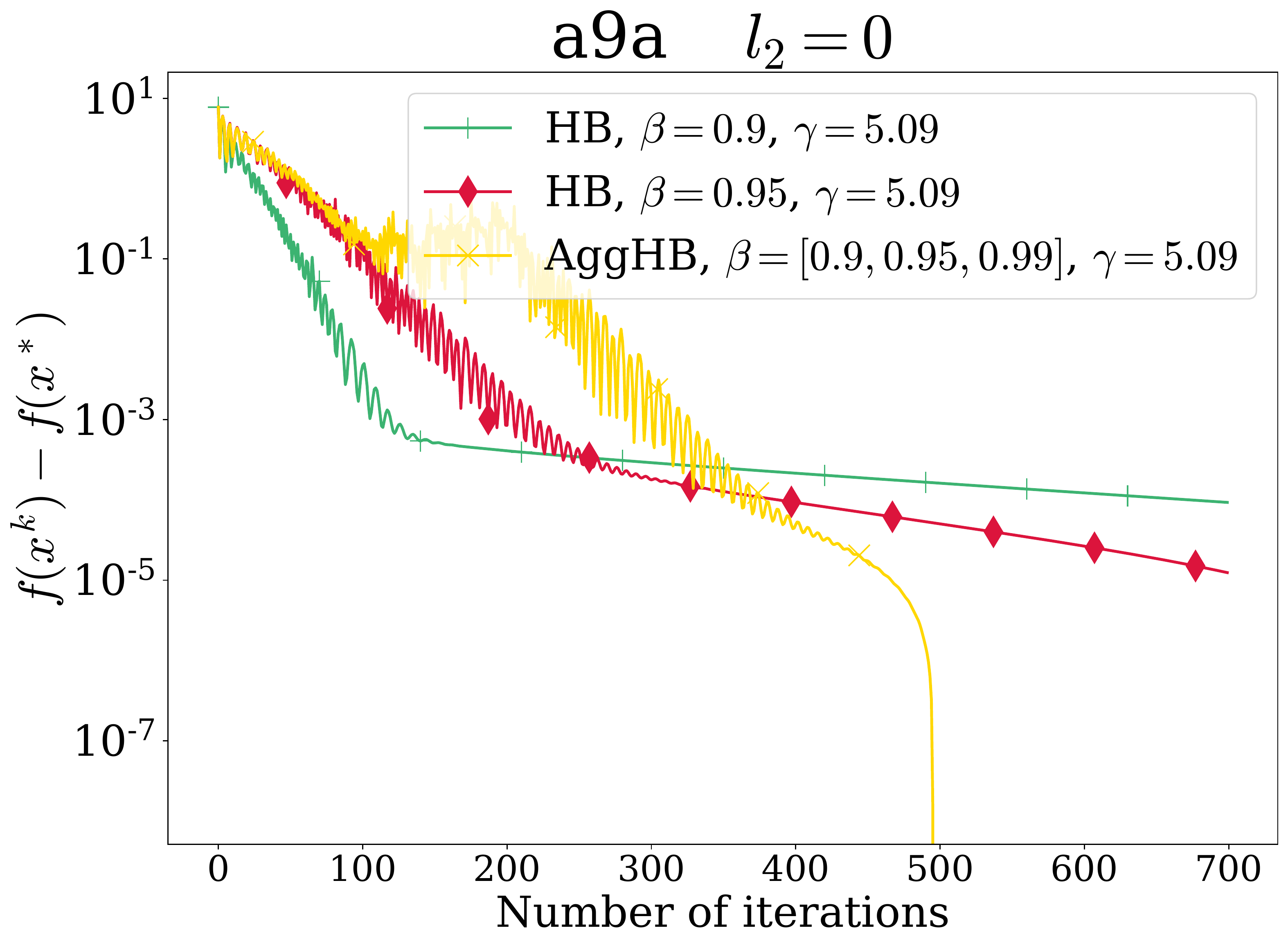}
    \includegraphics[width=0.325\textwidth]{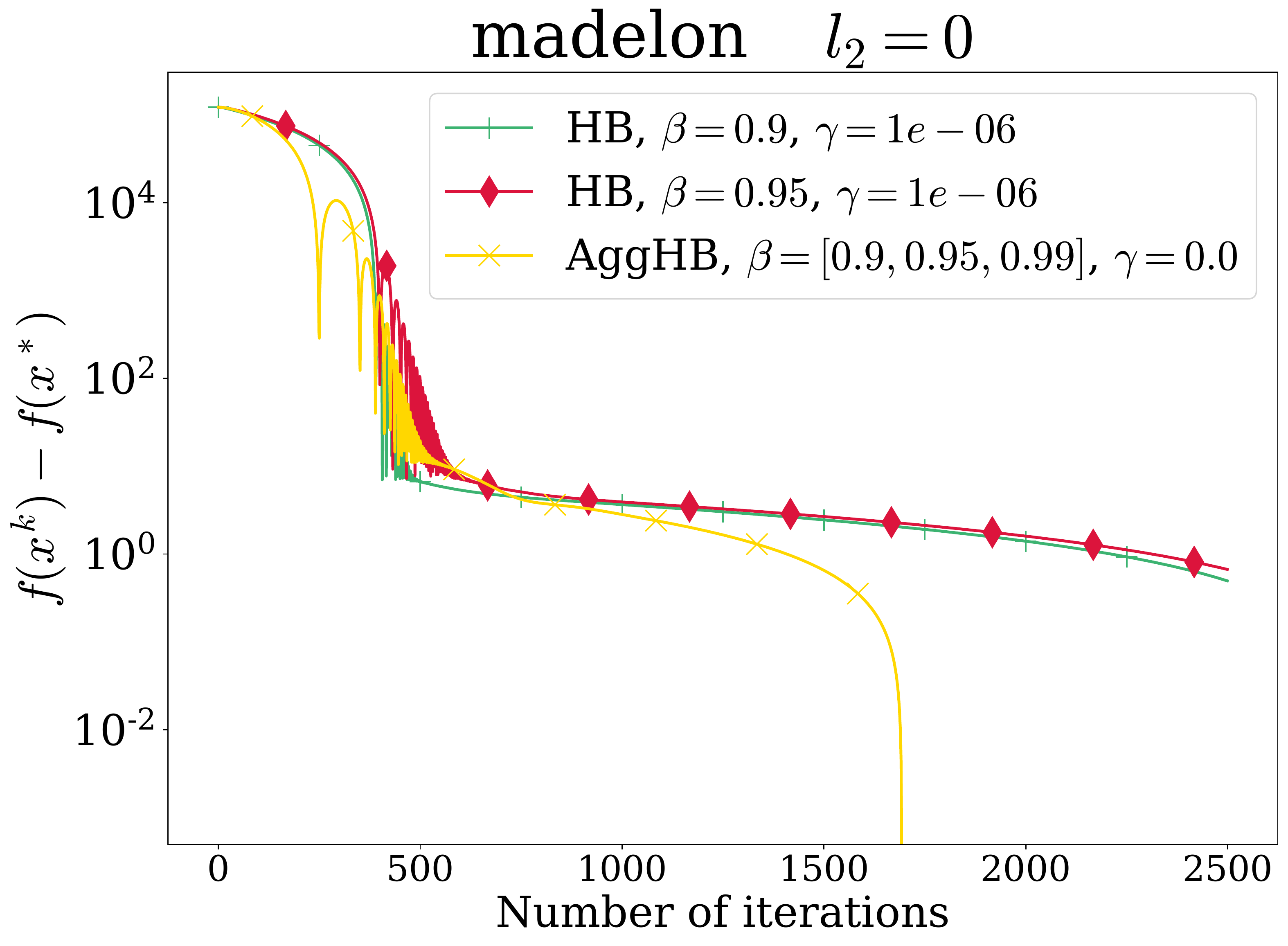}
    \includegraphics[width=0.325\textwidth]{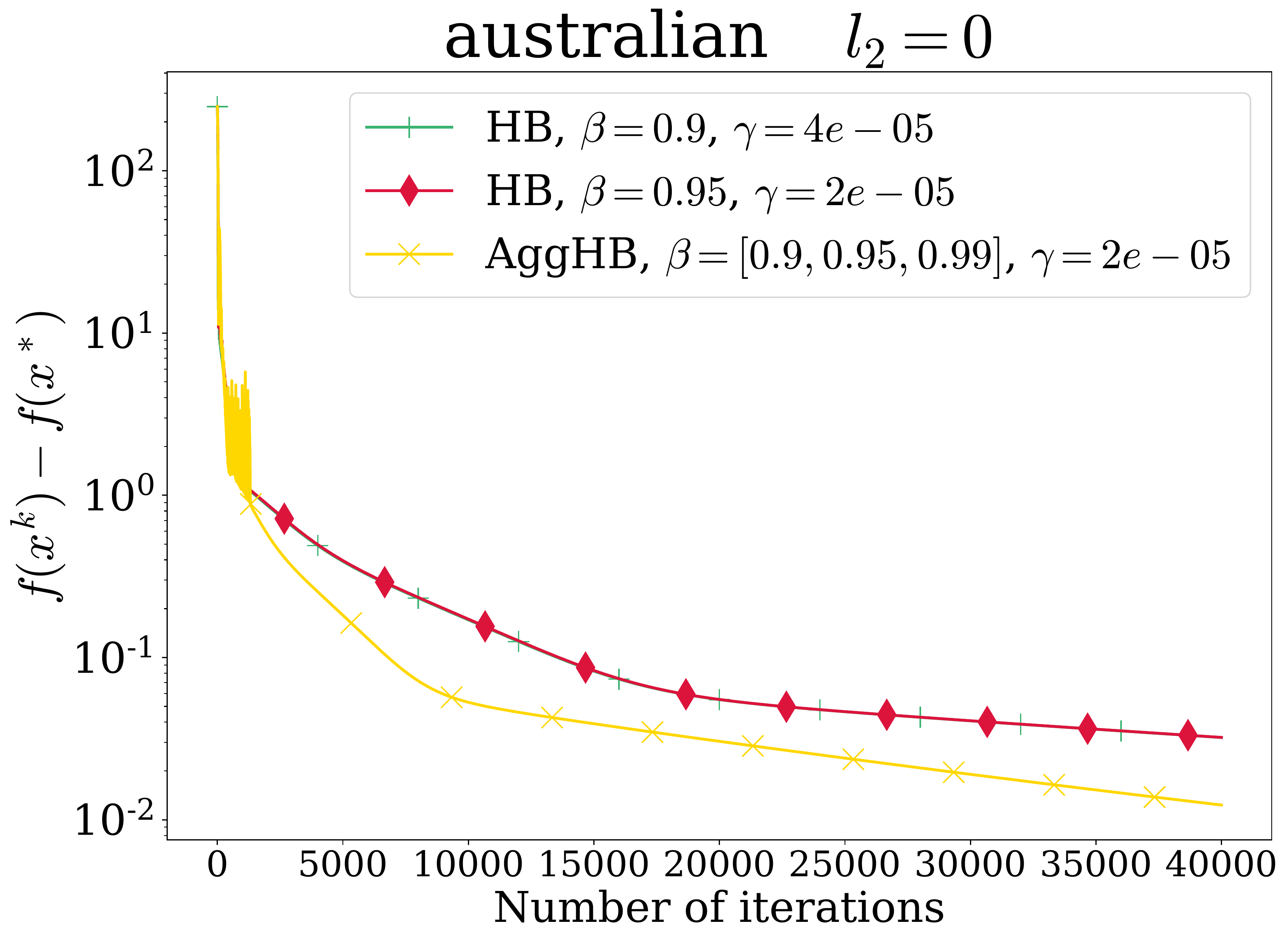}
    \includegraphics[width=0.325\textwidth]{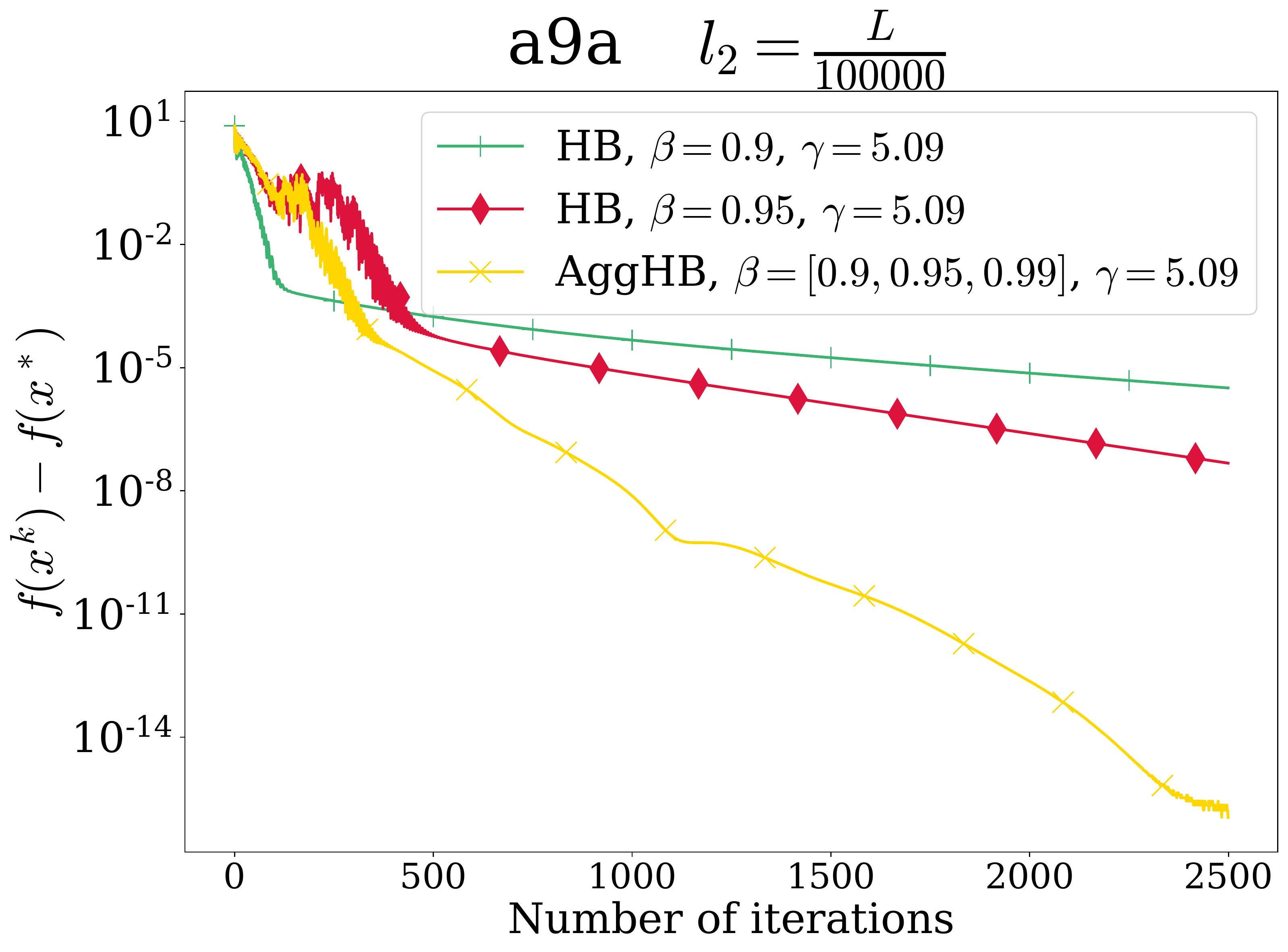}
    \includegraphics[width=0.325\textwidth]{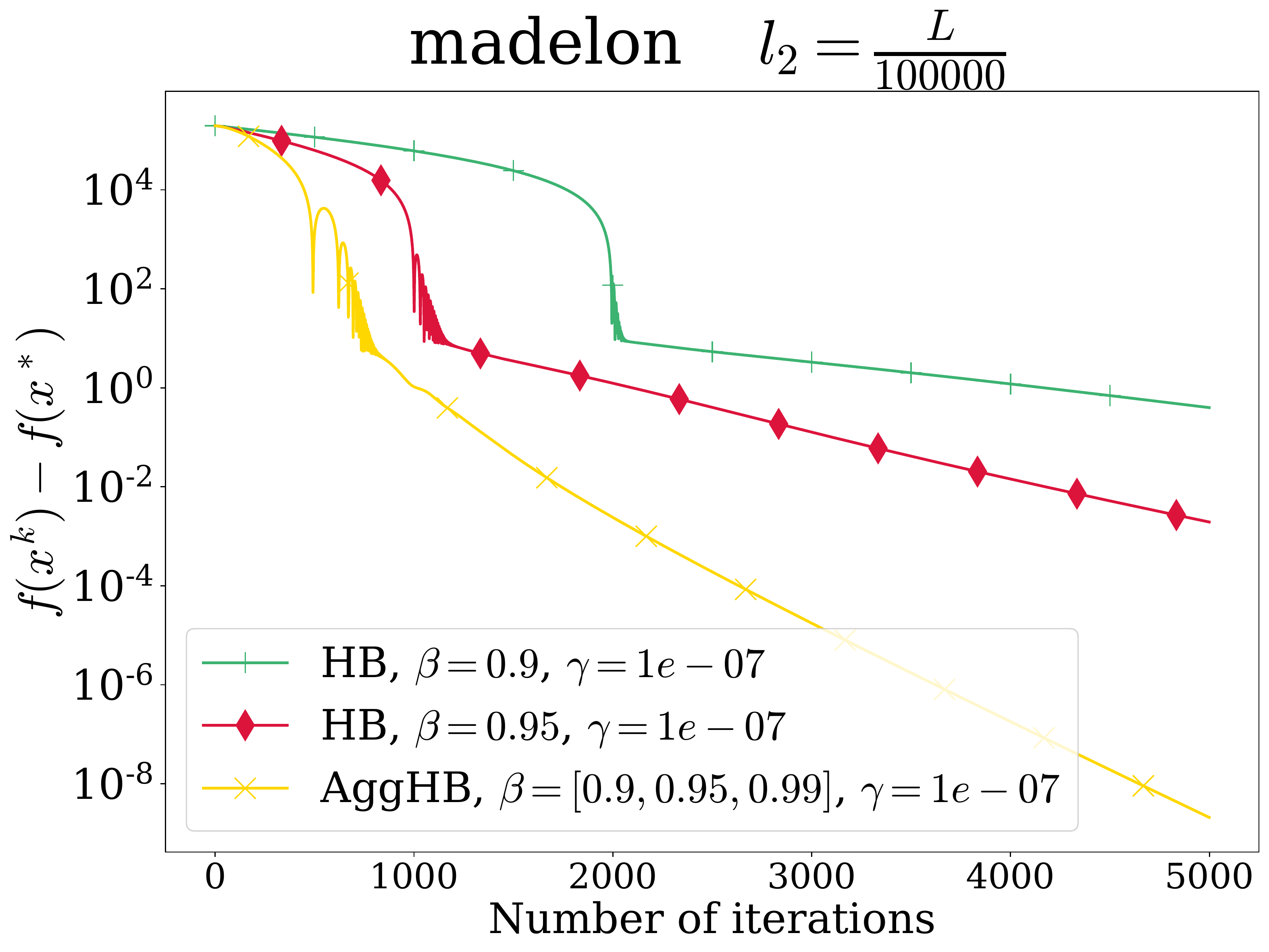}
    \includegraphics[width=0.325\textwidth]{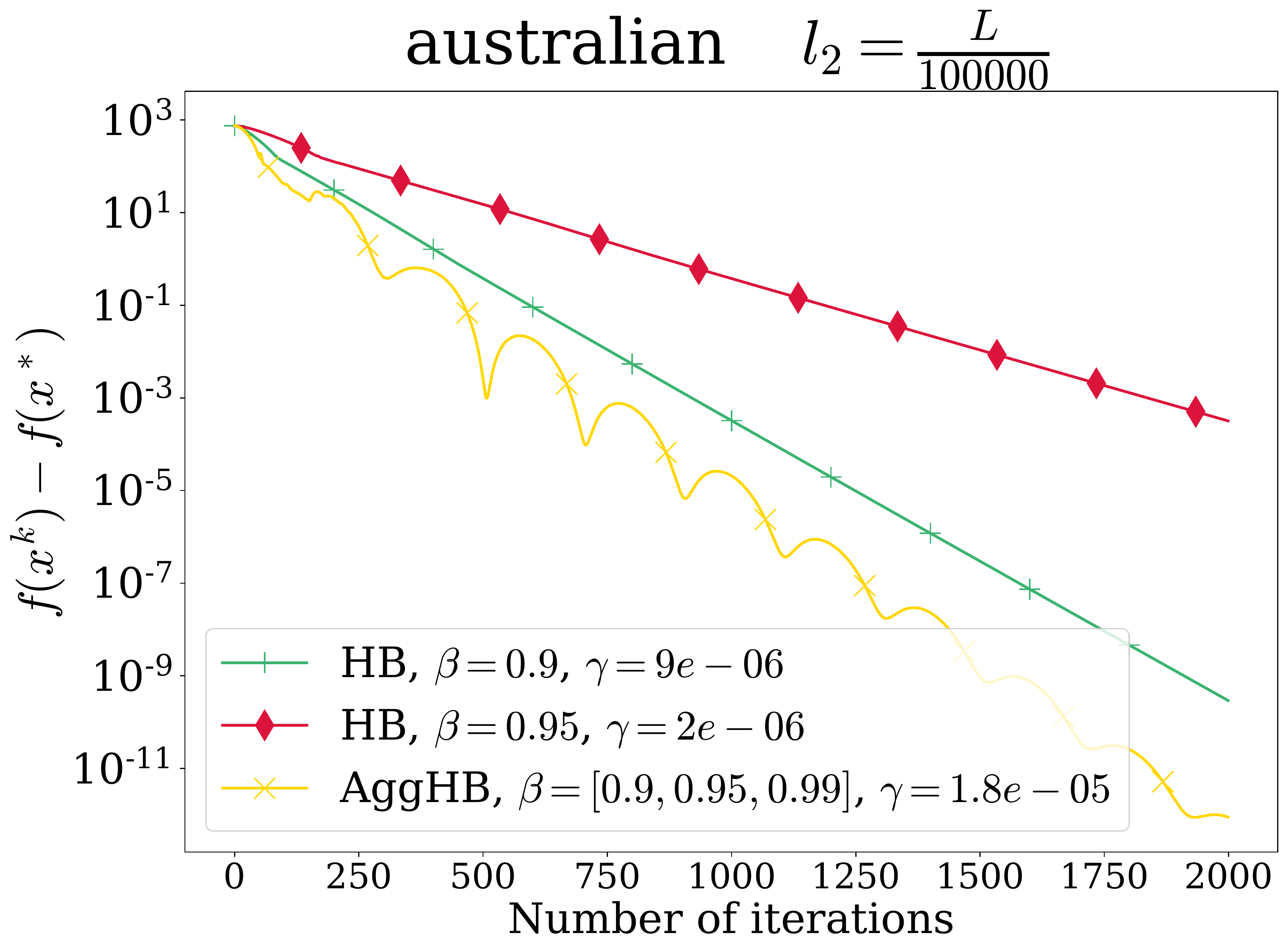}
    \includegraphics[width=0.325\textwidth]{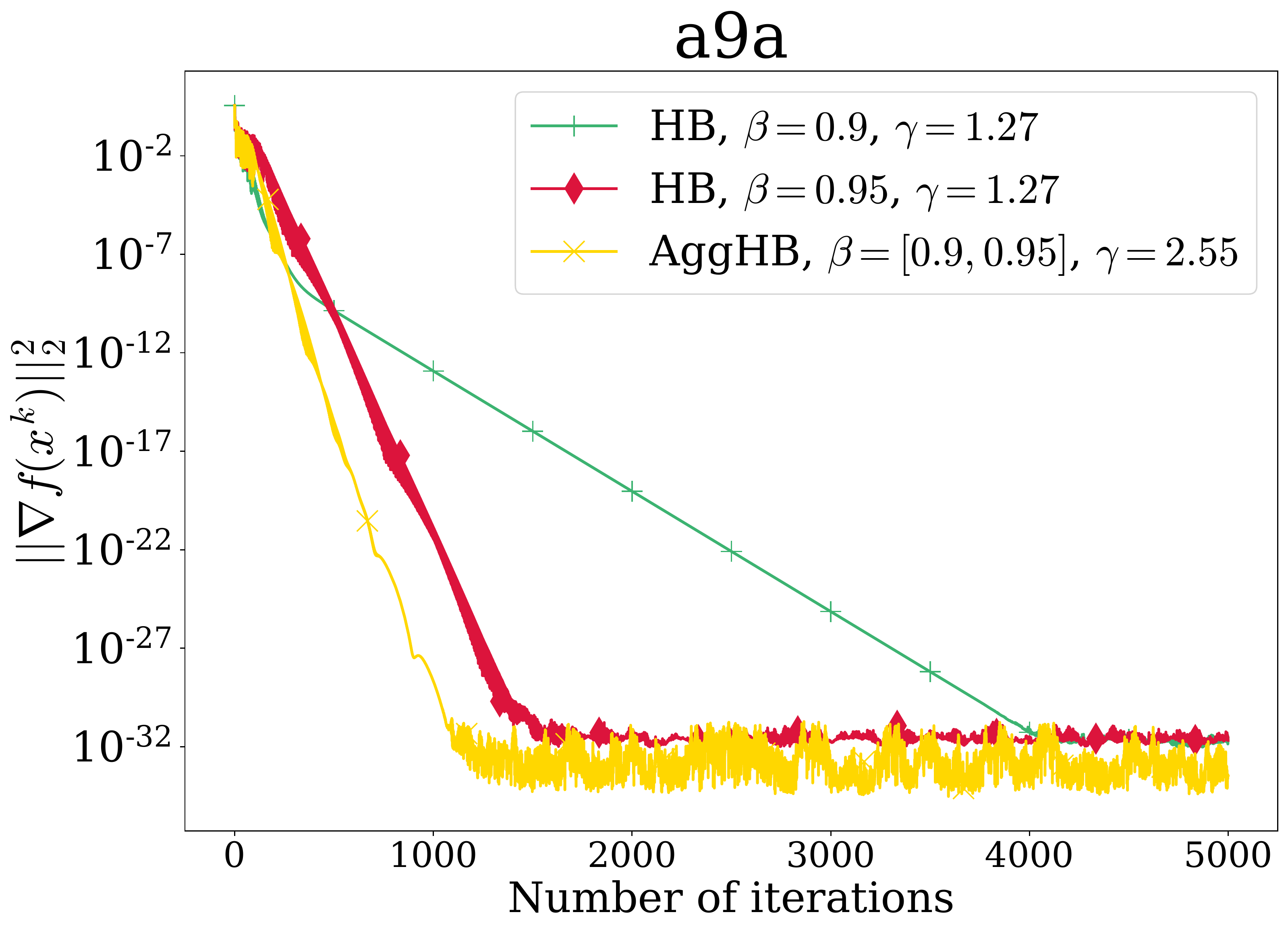}
    \includegraphics[width=0.325\textwidth]{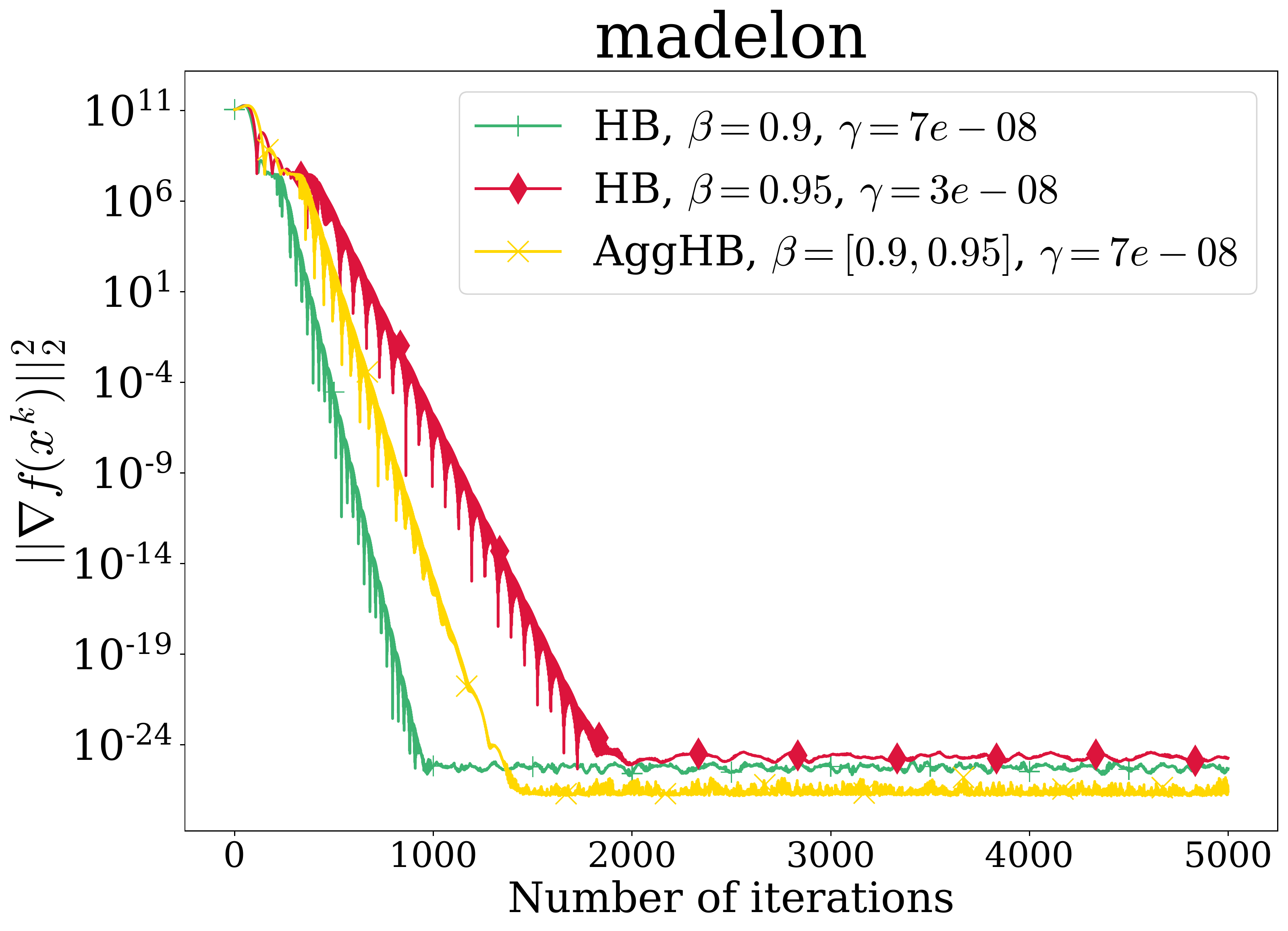}
    \includegraphics[width=0.325\textwidth]{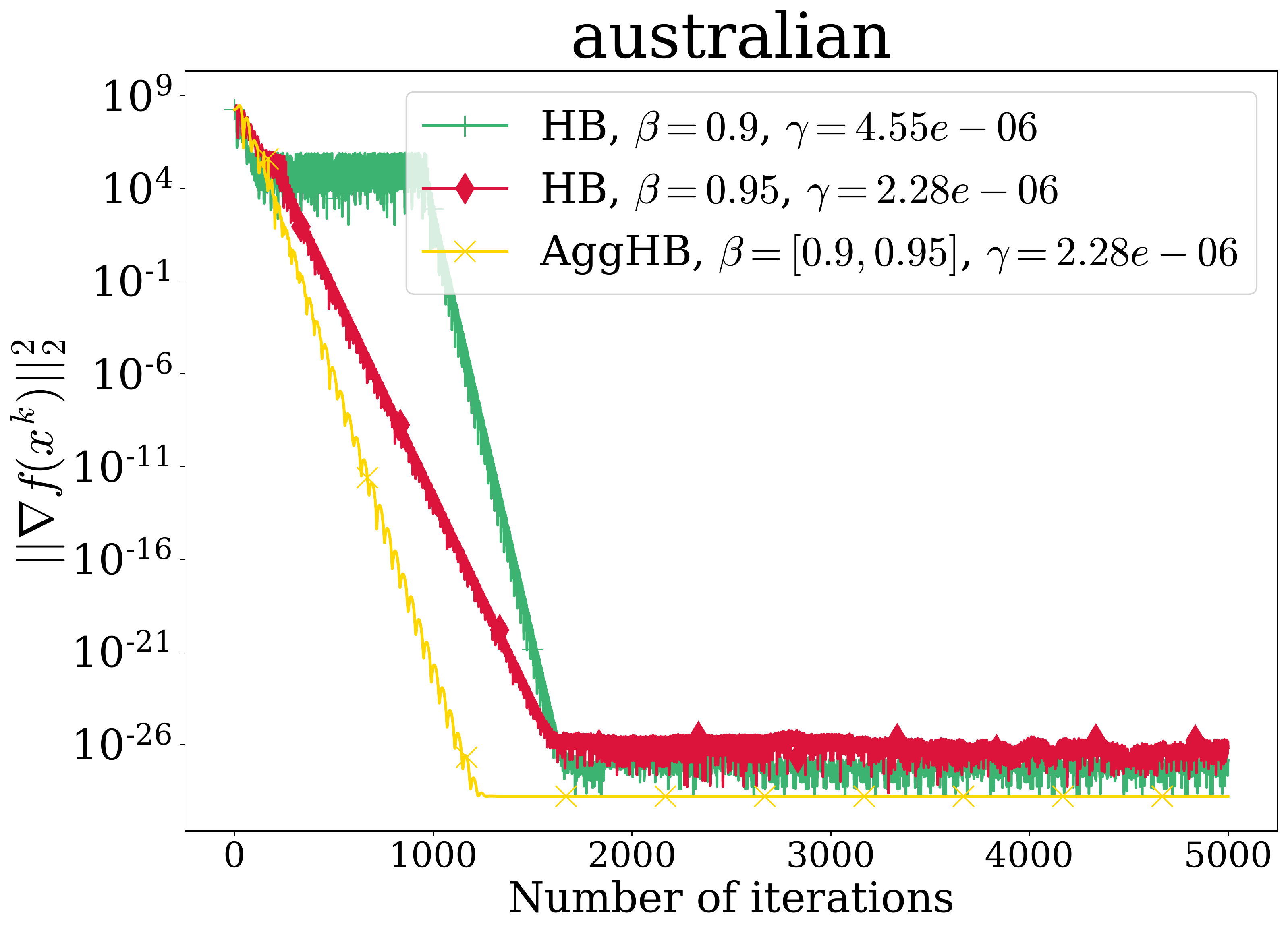}
\caption{Trajectories of \algname{HB} and \algname{AHB} with different momentum parameters $\beta$ applied to solve logistic regression problem with $\ell_2$-regularization (the first two rows) and non-convex regularization (the third row) for {\tt a9a}, {\tt madelon}, and {\tt australian} datasets. Stepsize $\gamma$ was tuned for each method.}
    \label{fig:cvx_strcvx}
\end{figure}

We compare the behavior of \algname{HB} and \algname{AggHB} on solving logistic regression problem with $\ell_2$-regularization and with special non-convex regularization:
\begin{eqnarray}
    &\min\limits_{x\in \R^n}&\left\{f(x) = \frac{1}{M}\sum\limits_{i=1}^M \log\left(1 + \exp\left(-y_i\cdot [Ax]_i\right)\right) + \frac{l_2}{2}\|x\|_2^2\right\}, \label{eq:logreg_l2}\\
    &\min\limits_{x\in \R^n}&\left\{f(x) = \frac{1}{M}\sum\limits_{i=1}^M \log\left(1 + \exp\left(-y_i\cdot [Ax]_i\right)\right) + \lambda \sum\limits_{j=1}^n\frac{x_j^2}{1+x_j^2}\right\}, \label{eq:logreg_noncvx}
\end{eqnarray}
where $M$ denotes the number of samples in the dataset, $A \in \R^{M\times n}$ is a ``feature matrix'', $y_1,\ldots,y_M \in \{-1,1\}$ are labels, and $l_2,\lambda \geq 0$ are the regularization parameters. One can show that $f(x)$ is $L$-smooth and $\mu$-strongly convex with $L = \tfrac{1}{4M}\lambda_{\max}(A^\top A) + l_2$ and $\mu = l_2$ in the first case, and $L$-smooth and non-convex with $L = \tfrac{1}{4M}\lambda_{\max}(A^\top A) + 2\lambda$. To construct the problems we use the following datasets from LIBSVM \cite{chang2011libsvm}: {\tt a9a} ($M = 32561$, $n = 123$), {\tt madelon} ($M = 2000$, $n = 500$), and {\tt australian} ($M = 690$, $n = 14$). Regularization parameter $l_2$ is either $0$ (convex problem) or $\tfrac{L}{100000}$ (strongly convex problem) and $\lambda$ is chosen as $\lambda = \tfrac{L}{1000}$. We run \algname{HB} with standard momentum parameters $\beta = 0.9, 0.95$ for both problems. \algname{AggHB} was tested with $m = 3$, $\beta_1 = 0.9$, $\beta_2 = 0.95$, $\beta_3 = 0.99$, and $\gamma_1 = \gamma_2 = \gamma_3 = \gamma$ for $\ell_2$-regularized problem and with $m = 2$, $\beta_1 = 0.9$, $\beta_2 = 0.95$, and $\gamma_1 = \gamma_2 = \gamma$. For each method we tune stepsize parameter $\gamma$ as follows: we choose $\gamma = \tfrac{a}{L}$ with the best $a \in \{2^{-6}, 2^{-5}, 2^{-4}, \ldots, 2^{8}\}$, i.e., the method achieves the best accuracy with the chosen $a$ from the considered set.

The results are shown in the Figure~\ref{fig:cvx_strcvx}. We observe that \algname{AggHB} outperforms \algname{HB} in all cases. In particular, for $\ell_2$-regularized problem the large value of $\beta_3$ does not slow down the convergence of \algname{AggHB}. In contrast, we observed that \algname{HB} performs relatively bad with $\beta = 0.99$. Next, in the experiments with non-convex regularization, \algname{AggHB} takes the best from two choices of momentum parameters.

\section{Conclusion}\label{sec:conclusion}
In this paper, we obtain the first convergence guarantees for \algname{AggHB} without assuming that the gradients of the objective function are uniformly bounded. In the special case when $m=1$, our results recover the known ones for \algname{HB} and outperform the corresponding guarantees for \algname{HB} with $\beta = \max_{i=1,m}\beta_i$ when $m > 1$. Our numerical results show the superiority of \algname{AggHB} to \algname{HB}. Together with the results from \cite{lucas2019aggregated} they indicate high practical potential of \algname{AggHB}.

\bibliography{references}

\appendix

\section{Missing Proofs from Section~\ref{sec:convergence_guarantees}}

\subsection{Proof of Theorem~\ref{thm:AggHB_non_convex}}

From $L$-smoothness of $f$ we have
\begin{eqnarray}
    f(\widetilde{x}_{k+1}) \le f(\widetilde{x}_k) - A \left\langle \nabla f(\widetilde{x}_k), \nabla f(x_k) \right\rangle + \frac{L A^2}{2} \|\nabla f(x_k)\|^2_2, 
    \label{eq:L_smoothness_non_cvx}
\end{eqnarray}
where $A=\frac{1}{m}\sum\limits_{i=1}^m \frac{\beta_{i}}{1-\beta_{i}}\gamma_{i}$. Next, we estimate a second term $-A\left\langle \nabla f(\widetilde{x}_k), \nabla f(x_k) \right\rangle$ in the previous expression:
\begin{eqnarray}
    -A\left\langle \nabla f(\widetilde{x}_k), \nabla f(x_k) \right\rangle &=& A \frac12 \left(\|\nabla f(\widetilde{x}_k) - \nabla f(x_k)\|^2_2 - \|\nabla f(\widetilde{x}_k)\|^2_2 - \|\nabla f(x_k)\|^2_2\right)  \notag\\
   &\overset{\eqref{eq:L_smoothness_non_cvx}}{\le}&  \frac{A}{2}(L^2 \|\widetilde{x}_k-x_k\|^2_2 - \|\nabla f(x_k)\|^2_2) \notag\\
  &\overset{\eqref{eq:virtual_iterates_AggHB}}{=}& \frac{AL^2}{2m^2} \left\| \sum\limits_{i=1}^m \frac{\beta_{i} \gamma_{i}}{1-\beta_{i}} V_{k-1}^{(i)}\right\|^2_2 - \frac{A}{2}\|\nabla f(x_k)\|^2_2.
  \label{eq:non_cvx_technical_1}
\end{eqnarray}
From \algname{AggHB} update rule we know that $V_k^{(i)}$ is linear combination of gradients: $V_k^{(i)} = \sum_{l=0}^k (\beta_{i})^l\nabla f(x_{k-l})$. Applying this to \eqref{eq:non_cvx_technical_1} we have
\begin{eqnarray}
    \frac{AL^2}{2m^2}\left\| \sum\limits_{i=1}^m \frac{\beta_{i} \gamma_{i}}{1-\beta_{i}} V_{k-1}^{(i)} \right\|^2_2 \le \frac{AL^2B}{2m^2} \sum\limits_{l=0}^{k-1}\left\| \nabla f (x_{k-1-l})\right\|^2_2 \sum\limits_{i=1}^m \frac{(\beta_{i})^l \gamma_{i}}{1-\beta_{i}},
    \label{eq:non_cvx_technical_2}
\end{eqnarray}
where $B = \sum_{l=0}^{k-1} \sum_{i=1}^m \tfrac{(\beta_{i})^l \gamma_{i}}{1-\beta_{i}} \le \sum_{i=1}^m \tfrac{\gamma_{i}}{(1-\beta_{i})^2}.$ Combining \eqref{eq:non_cvx_technical_1}, \eqref{eq:non_cvx_technical_2}, we continue the derivation from \eqref{eq:L_smoothness_non_cvx}:
\begin{eqnarray}
    f(\widetilde{x}_{k+1}) &\le& f(\widetilde{x}_k) - \frac{A}{2} (1-LA) \|\nabla f(x_k)\|^2_2 + \frac{AL^2B}{2m^2} \sum\limits_{l=0}^{k-1} \|\nabla f(x_{k-1-l})\|^2_2 \sum\limits_{i=1}^m \frac{(\beta_{i})^l \gamma_{i}}{1-\beta_{i}} \notag\\
    &\le&  f(\widetilde{x}_k) - \frac{A}{2}(1-LA)\|\nabla f(x_k)\|^2_2 \notag\\
    &&  + \frac{AL^2}{2m^2} \left(\sum\limits_{i=1}^m \frac{\gamma_{i}}{(1-\beta_{i})^2} \right) \sum\limits_{l=0}^{k-1}\|\nabla f(x_l)\|^2_2 \sum\limits_{i=1}^m \frac{(\beta_{i})^{k-1-l}\gamma_{i}}{1-\beta_{i}} \notag\\ &\le& f(\widetilde{x}_k) -\frac{A}{2}(1-LA)\|\nabla f(x_k)\|^2_2 \notag \\
    & & + \frac{AL^2}{2m^2}\!\left(\sum\limits_{i=1}^m \frac{\gamma_{i}}{(1-\beta_{i})^2}\right)\! \left(\max\limits_{i=1,m} \frac{\gamma_{i}}{1-\beta_{i}} \right)\! \sum\limits_{l=0}^{k-1}\sum\limits_{i=1}^m (\beta_{i})^{k-1-l} \|\nabla f(x_l)\|^2_2.
    \label{eq:non_cvx_sum}
\end{eqnarray}
Summing up \eqref{eq:non_cvx_sum} for $k=0,1,\ldots, K$ we get

\begin{eqnarray*}
    f(\widetilde{x}_{k+1}) &\le& f(\widetilde{x}_0) +
    \sum\limits_{k=1}^K  \left(\frac{LA^2-A}{2}\right)\|\nabla f(x_{k})\|^2_2 \\
    && +  \sum\limits_{k=1}^K\! \left( \frac{AL^2}{2m^2}\! \left( \sum\limits_{i=1}^m \frac{\gamma_{i}}{(1\!-\!\beta_i)^2}\right)\! \left(\max\limits_{i=1,m} \frac{\gamma_{i}}{1\!-\!\beta_i} \right)\! \sum\limits_{i=1}^m\! \sum\limits_{l=k+1}^{K-1} \beta_i^{l-1-k} \right)\! \|\nabla f(x_{k})\|^2_2 \\
    &\le& f(\widetilde{x}_0) + \sum\limits_{k=1}^K  \left(\frac{LA^2-A}{2}\right)\|\nabla f(x_{k})\|^2_2 \\
    && +  \sum\limits_{k=1}^K \left( \frac{AL^2}{2m^2} \left( \sum\limits_{i=1}^m \frac{\gamma_{i}}{(1-\beta_i)^2}\right) \left(\max\limits_{i=1,m} \frac{\gamma_{i}}{1-\beta_i} \right) \sum\limits_{i=1}^m \frac{1}{1-\beta_i} \right) \|\nabla f(x_{k})\|^2_2 \\
    &=& f(\widetilde{x}_0) + \sum\limits_{k=1}^K \left( \left(\frac{LA^2-A}{2}\right)
    - \frac{ACDEL^2}{2m^2} \right) \|\nabla f(x_{k})\|^2_2, 
\end{eqnarray*}
where
$A = \frac{1}{m}\sum\limits_{i=1}^m\frac{\beta_i \gamma_{i}}{1-\beta_i}, \; C = \sum\limits_{i=1}^m \frac{\gamma_{i}}{(1-\beta_i)^2}, \; D = \max\limits_{i=1,m} \frac{\gamma_{i}}{1-\beta_i}, \; E = \sum\limits_{i=1}^m \frac{1}{1-\beta_i}$. Finally, by choosing sufficiently small $ \gamma_i $ one can ensure that $-\frac{A}{2}\left( 1 - \frac{CDEL^2}{2m^2} - LA\right) \le 0$ and get \eqref{eq:main_result_non_cvx}

\subsection{Proof of
Corollary~\ref{cor:AggHB_complexity_non_cvx}}

From Theorem~\ref{thm:AggHB_non_convex} we have
$\min\limits_{k=1,K} \|\nabla f(x_{k})\|^2_2 \le \frac{2}{K} \frac{f(x_0)-f_{\inf}}{A\left(1 - \frac{CDEL^2}{m^2} - LA\right)}$ for all $K \ge 1.$
This upper bound implies that to achieve $\min\limits_{k=1,K}\|\nabla f(x_k)\|^2_2 \le \varepsilon^2$, the method requires
    \begin{eqnarray}
        K = \cO\left( \frac{\Delta_0}{A\left(1 - \frac{CDEL^2}{m^2} - LA\right)\varepsilon^2}\right)\quad \text{iterations},
        \label{eq:cor_non_cvx_tech_1}
    \end{eqnarray}
where $\Delta_0 = f(x_0)-f_{\inf}.$ It remains to estimate the denominator in the above complexity bound. To do that, we introduce new constants $\widetilde{\beta}$ and $\hat{\beta}$ satisfying the following conditions
    \begin{equation*}
        \frac{1}{m}\sum\limits_{i=1}^m \frac{\beta_i}{(1-\beta_i)^2} = \frac{\widetilde{\beta}}{(1-\widetilde{\beta})^2}, \quad \frac{1}{m}\sum\limits_{i=1}^m \frac{1}{1-\beta_i} = \frac{1}{1-\hat{\beta}}.
    \end{equation*}
Next, using the formulas for $\widetilde{\beta}$, $\hat{\beta}$ and assuming $\gamma_i=\gamma$, we get new expressions for constants $A, C, D, E$:
    \begin{eqnarray*}
        A &=& \frac{1}{m}\sum\limits_{i=1}^m\frac{\beta_i \gamma_{i}}{1-\beta_i} = \frac{\gamma \hat{\beta}}{1-\hat{\beta}},\quad C = \sum\limits_{i=1}^m \frac{\gamma_{i}}{(1-\beta_i)^2}=\gamma m \left(\frac{\widetilde{\beta}}{(1-\widetilde{\beta})^2} + \frac{1}{1-\hat{\beta}}\right), \notag\\
        D &=& \max\limits_{i=1,m} \frac{\gamma_{i}}{1-\beta_i} = \frac{\gamma}{1-\max\limits_{i=1,m} \beta_i},\quad E = \sum\limits_{i=1}^m \frac{1}{1-\beta_i} = \frac{m}{1-\hat{\beta}}.
    \end{eqnarray*}
Then, condition \eqref{eq:condition_non_cvx}, which is equivalent to 
$1 - \tfrac{CDEL^2}{m^2} - LA > 0$, can be written as
    \begin{eqnarray*}
        1 - \gamma^2L^2\left(\frac{\widetilde{\beta}}{(1-\widetilde{\beta})^2} + \frac{1}{1-\hat{\beta}}\right)\cdot\frac{1}{\left(1-\max\limits_{i=1,m}\beta_i\right)(1-\hat{\beta})} - \gamma L\frac{\hat{\beta}}{1-\hat{\beta}} > 0.
    \end{eqnarray*}
To derive the complexity stated in the corollary, we choose $\gamma$ such that
    \begin{eqnarray*}
        \frac{1}{2} - \gamma^2L^2\left(\frac{\widetilde{\beta}}{(1-\widetilde{\beta})^2} + \frac{1}{1-\hat{\beta}}\right)\cdot\frac{1}{\left(1-\max\limits_{i=1,m}\beta_i\right)(1-\hat{\beta})} - \gamma L\frac{\hat{\beta}}{1-\hat{\beta}} \geq 0,
    \end{eqnarray*}
implying \eqref{eq:condition_non_cvx}. One can show (see Lemma 5 from \cite{richtarik2021ef21}) that this condition holds for
    \begin{equation*}
        \gamma = \frac{1}{L\left(\frac{2\hat{\beta}}{1-\hat{\beta}} + \sqrt{2\left(\frac{\widetilde{\beta}}{(1-\widetilde{\beta})^2} + \frac{1}{1-\hat{\beta}}\right)\frac{1}{\left(1-\max\limits_{i=1,m}\beta_i\right)(1-\hat{\beta})}}\right)}
    \end{equation*}
Plugging this value of $\gamma$ in \eqref{eq:cor_non_cvx_tech_1} and using $1 - \tfrac{CDEL^2}{m^2} - LA \geq \tfrac{1}{2}$, we finally obtain
    \begin{eqnarray*}
        K = \cO\left(\frac{L\Delta_0}{\varepsilon^2} + \frac{L\Delta_0 \sqrt{\left(\frac{\widetilde{\beta}(1-\hat\beta)}{(1-\widetilde{\beta})^2} + 1\right)\frac{1}{\left(1-\max\limits_{i=1,m}\beta_i\right)\hat\beta^2}}}{\varepsilon^2}\right).
    \end{eqnarray*}

\subsection{Proof of Lemma~\ref{lem:one_iter_progress_AggHB}}

Applying the virtual iterates determined in \eqref{eq:virtual_iterates_AggHB}, we obtain 
    \begin{eqnarray}
        \|\widetilde{x}_{k+1} - x_*\|^2_2 &=& \|\widetilde{x}_k - x_*\|^2_2 - 2F\langle \widetilde{x}_k - x_*, \nabla f(x_k) \rangle + F^2\|\nabla f(x_k)\|^2_2 \notag\\
        &=& \|\widetilde{x}_k - x_*\|^2_2 - 2F\langle x_k - x_*, \nabla f(x_k) \rangle - 2F\langle \widetilde{x}_k - x_k, \nabla f(x_k) \rangle\notag\\
        &&\quad + F^2\|\nabla f(x_k)\|^2_2 \label{eq:one_iter_AggHB_technical_1}.
    \end{eqnarray}
    From $\mu$-strong convexity and $L$-smoothness of $f$ we have (e.g., see \cite{nesterov2018lectures})
    \begin{eqnarray}
        \langle x_k - x_*, \nabla f(x_k) \rangle &\ge& f(x_k) - f(x_*) + \frac{\mu}{2}\|x_k - x_*\|^2 \notag\\
        \|\nabla f(x_k)\|^2 &\le& 2L\left(f(x_k) - f(x_*)\right). \label{eq:L_smoothness_cor}
    \end{eqnarray}
    Using these inequalities for \eqref{eq:one_iter_AggHB_technical_1} we get
    \begin{eqnarray*}
        \|\widetilde{x}_{k+1} - x_*\|^2_2 &\le& \|\widetilde{x}_k - x_*\|^2_2 - \mu F \|x_k - x_*\|^2_2 - 2F\left(f(x_k) - f(x_*)\right) \\
        &&\quad - 2F \langle \widetilde{x}_{k} - x_k, \nabla f(x_k)  \rangle + F^2\|\nabla f(x_k)\|^2_2.
    \end{eqnarray*}
    Firstly, we evaluate the second term $-\mu F \|x_k - x_*\|^2_2$ using that $\|a+b\|_2^2 \le 2\|a\|_2^2 + 2\|b\|_2^2$ for all $a,b \in \R^n$ as follows
    \begin{equation*}
        - \mu F \|x_k - x_*\|^2_2 \le -\frac{\mu F}{2}\|\widetilde{x}_k - x_*\|^2_2 +\mu F\|x_k - \widetilde{x}_k\|^2_2.
    \end{equation*}
    Secondly, we estimate the fourth term $ - 2F \langle \widetilde{x}_{k} - x_k, \nabla f(x_k) $ using Fenchel-Young inequality\footnote{$|\la a, b\ra| \le \frac{\|a\|_2^2}{2\lambda} + \frac{\lambda\|b\|_2^2}{2}$ for all $a,b\in\R^n$ and $\lambda > 0$.} and get
    \begin{eqnarray*}
        - 2F \langle \widetilde{x}_{k} - x_k, \nabla f(x_k) &\le& -2LF \|\widetilde{x}_{k} - x_k\|^2_2 + \frac{F}{2L}\|\nabla f(x_k)\|^2_2\\
        &\overset{\eqref{eq:L_smoothness_cor}}{\le}& -2LF \|\widetilde{x}_{k} - x_k\|^2_2 + F\|\left(f(x_k) - f(x_*)\right).
    \end{eqnarray*}
    Combining the results above, we finish the proof
    \begin{eqnarray*}
        \|\widetilde{x}_{k+1} - x_*\|^2
        &\overset{\eqref{eq:virtual_iterates_AggHB},\eqref{eq:one_iter_progress_AggHB}}{\le}& \left(1 - \frac{\mu F}{2}\right)\|\widetilde{x}_k - x_*\|^2_2 - \frac{F}{2}\left(f(x_k) - f(x_*)\right) + 3LF\|x_k - \widetilde{x}_k\|^2_2.
    \end{eqnarray*}

\subsection{Proof of Lemma~\ref{lem:weighted_sum_of_momentums}}

From \algname{AggHB} update rule we know that $V_k^{(i)}$ is linear combination of gradients: $V_k^{(i)} = \sum_{t=0}^k \beta_{i}^t\nabla f(x_{k-t})$. Next, by the definition of $\widetilde{x}_k$ we have
\begin{eqnarray}
       \|x_{k+1} - \widetilde{x}_{k+1}\|^2_2 &=& \left\|\frac{1}{m}\sum\limits_{i=1}^m \frac{\beta_i \gamma_i}{1-\beta_i} V_k^{(i)}\right\|^2_2 
       =\left\|\frac{1}{m}\sum\limits_{i=1}^m \frac{\beta_i \gamma_i}{1-\beta_i}\sum\limits_{t=0}^k \beta_i^t \nabla f(x_{k-t})\right\|^2_2 \notag\\
       &=& \left\|\sum\limits_{t=0}^k \left( \frac{1}{m} \sum\limits_{i=1}^m \frac{\beta_i^{t+1} \gamma_i}{1-\beta_i}\right)\nabla f(x_{k-t})\right\|^2_2 \notag\\
       &=& \left\|\sum\limits_{t=0}^k \left( \frac{1}{m} \sum\limits_{i=1}^m \frac{\beta_i^{k-t+1} \gamma_i}{1-\beta_i}\right)\nabla f(x_{t})\right\|^2_2.
       \label{eq:technical_1}
\end{eqnarray}
Define constant $B_k$ as following
    \begin{eqnarray*}
        B_k = \sum\limits_{t=0}^k \frac{1}{m} \sum\limits_{i=1}^m \frac{\beta_i^{k-t+1}\gamma_i}{1-\beta_i} = \frac{1}{m} \sum\limits_{i=1}^m \frac{\beta_i \gamma_i}{1-\beta_i}\sum\limits_{t=0}^{k}\beta_i^{k-t} = \frac{1}{m}\sum\limits_{i=1}^m \frac{\beta_i \gamma_i (1-\beta_i^{k+1})}{(1-\beta_i)^2}.
    \end{eqnarray*}
Using this, we continue  the derivation from \eqref{eq:technical_1}
    \begin{eqnarray}
       \|x_{k+1} - \widetilde{x}_{k+1}\|^2_2 &=& 
       B_k^2 \cdot \left\|\sum\limits_{t=0}^k \frac{1}{B}\left( \frac{1}{m}\sum\limits_{i=1}^m \frac{\beta_i^{k-t+1}\gamma_i}{1-\beta_i}\right)\nabla f(x_t)\right\|^2_2 \notag\\
       &\overset{\text{Jensen's inequality}}{\le}&
       B_k^2 \cdot \sum\limits_{t=0}^k \frac{1}{B} \left( \frac{1}{m} \sum\limits_{i=1}^m \frac{\beta_i^{k-t+1}\gamma_i}{1-\beta_i}\right)\|\nabla f(x_t)\|^2_2  \notag\\
       &=& B_k \cdot \sum\limits_{t=0}^k \left( \frac{1}{m} \sum\limits_{i=1}^m \frac{\beta_i^{k-t+1}\gamma_i}{1-\beta_i}\right)\|\nabla f (x_t)\|^2_2 \notag\\
       &\overset{\eqref{eq:constants_AggHB_2}, \; B_k \le B_{k+1}}{\le}&
       B_K F \cdot \sum\limits_{t=0}^k \max\limits_{i=1,m}\beta_i^{k-t+1}\|\nabla f(x_t)\|^2_2.
       \label{eq:technical_2}
    \end{eqnarray}
For simplicity, we denote $B_K \equiv B$.
Summing up these inequalities for $k=0,1,\ldots, K$ with weights $w_k = \left(1 - \frac{\mu F}{2}\right)^{-(k+1)}$, we get
    \begin{eqnarray}
        3LF \sum\limits_{k=0}^K w_k \|x_k-\widetilde{x}_k\|^2_2 &\le& 3LBF^2 \cdot \sum\limits_{k=0}^K \sum\limits_{t=0}^{k-1}w_k \max\limits_{i=1,m} \beta_i^{k-t}\|\nabla f(x_t)\|^2_2 \notag\\
        &\le& 3LBF^2 \cdot \sum\limits_{k=0}^K \sum\limits_{t=0}^{k}w_k \max\limits_{i=1,m}\beta_i^{k-t}\|\nabla f(x_t)\|^2_2. \label{eq:weighted_sum_momentums_technical_1}
    \end{eqnarray}
    Next, we estimate $w_k$ using that $(1 - \nicefrac{q}{2})^{-1} \le 1+q$ for any $q\in (0,1]$: for all $t = 0,1,\ldots,k$
    \begin{equation}
        w_k = \left(1\! -\! \frac{\mu F}{2}\right)^{-(k-t)}\!w_{t} \leq \left(1\! +\! \mu F\right)^{k-t}\!w_t \overset{\eqref{eq:params_AggHB_2}}{\le} \left(1 \!+
        \!\frac{1\!-\!\max\limits_{i=1,m}\beta_i}{2}\right)^{k-t}\!w_t. \notag
    \end{equation}
    Using an inequality above and $\left(1 + \nicefrac{q}{2}\right)(1-q) \le 1-\nicefrac{q}{2}$ for $q = 1-\max_{i=1,m}\beta_i$, we continue the previous derivation \eqref{eq:weighted_sum_momentums_technical_1}
    \begin{align}
        3LF\!\sum\limits_{k=0}^K\! w_k\! \|x_k\!-\!\widetilde{x}_k\|^2_2 &\le 3LBF^2\! \sum\limits_{k=0}^K\!\sum\limits_{t=0}^k\! w_t\! \|\nabla f(x_t)\|^2_2\! \left(\!1\! +\! \frac{1-\max\limits_{i=1,m}\beta_i}{2}\!\right)^{k-t}\! \max\limits_{i=1,m}\! \beta_i^{k-t} \notag\\
        &\le
        3LBF^2 \sum\limits_{k=0}^K\sum\limits_{t=0}^k w_t \|\nabla f(x_t)\|^2_2 \left(1\! -\! \frac{1\!-\!\max\limits_{i=1,m}\!\beta_i}{2}\right)^{k-t} \notag\\
        &\le  3LBF^2 \left(\sum\limits_{k=0}^{K}w_k \|\nabla f(x_k)\|^2_2\right) \left(\sum\limits_{k=0}^{\infty}\left(1\!-\!\frac{1\!-\!\max\limits_{i=1,m}\!\beta_i}{2}\right)^k\right) \notag\\
        &= \frac{6LBF^2}{1-\max\limits_{i=1,m}\beta_i} \sum\limits_{k=0}^K w_k \|\nabla f(x_k)\|^2_2 \notag\\
        &\le \frac{12L^2BF^2}{1-\max\limits_{i=1,m}\beta_i} \sum\limits_{k=0}^K w_k (f(x_k) - f(x_*)).
        \label{eq:weighted_sum_momentums_technical_2}
    \end{align}
    We take parameters $\gamma_i, \beta_i$ \eqref{eq:params_AggHB_2} implying \eqref{eq:constants_AggHB_2}. Combining this with the last result \eqref{eq:weighted_sum_momentums_technical_2},
    we obtain \eqref{eq:weighted_sum_of_momentums}.

\subsection{Proof of Theorem~\ref{thm:AggHB_main_result}}

Using Lemma~\ref{lem:one_iter_progress_AggHB} we get
    \begin{equation*}
        \frac{F}{2}\left(f(x_k) - f(x_*)\right) \le \left(1 - \frac{\mu F}{2}\right)\|\widetilde{x}_k - x_*\|_2^2 - \|\widetilde{x}_{k+1} - x_*\|_2^2 + 3LF\|x_k - \widetilde{x}_k\|_2^2.
    \end{equation*}
    Summing up these inequalities for $k = 0,1,\ldots,K$ with weights $w_k = \left(1 - \frac{\mu F}{2}\right)^{-(k+1)}$, we have
    \begin{eqnarray*}
        \frac{F}{2}\sum\limits_{k=0}^K w_k\left(f(x_k) - f(x_*)\right) &\le& \sum\limits_{k=0}^K\left(w_k\left(1 - \frac{\mu F}{2}\right)\|\widetilde{x}_k - x_*\|_2^2 - w_k\|\widetilde{x}_{k+1} - x_*\|_2^2\right)\\
        &&\quad + 3LF\sum\limits_{k=0}^K w_k\|x_k - \widetilde{x}_k\|_2^2\\
        &\overset{\eqref{eq:weighted_sum_of_momentums}}{\le}& \sum\limits_{k=0}^K\left(w_{k-1}\|\widetilde{x}_k - x_*\|_2^2 - w_k\|\widetilde{x}_{k+1} - x_*\|_2^2\right)\\
        &&\quad + \frac{F}{4}\sum\limits_{k=0}^K w_k\left(f(x_k) - f(x_*)\right)\\
        &\le& \|x_0 - x_*\|_2^2 + \frac{F}{4}\sum\limits_{k=0}^K w_k\left(f(x_k) - f(x_*)\right).
    \end{eqnarray*}
    Rearranging and multiplying by $\frac{1}{W_K} = \frac{1}{\sum_{k=0}^K w_k}$ this inequality, we have
    \begin{equation*}
        \frac{1}{W_K}\sum\limits_{k=0}^K w_k\left(f(x_k) - f(x_*)\right) \le \frac{4\|x_0 - x_*\|_2^2}{F W_K}.
    \end{equation*}
    Next, we obtain \eqref{eq:AggHB_main_result} by using Jensen's inequality:
    \begin{equation*}
        f(\overline{x}_K) \le \frac{1}{W_K}\sum\limits_{k=0}^K w_k f(x_k).
    \end{equation*}
     
    In strongly convex case ($\mu > 0$), we have $W_K \ge w_{K-1} = \left(1 - \frac{\mu F}{2}\right)^{-K}$, hence \eqref{eq:AggHB_main_result_str_cvx} holds. In convex case ($\mu = 0$), $W_K = {K+1} > K$ that implies \eqref{eq:AggHB_main_result_cvx}.
    
\subsection{Proof of
Corollary~\ref{cor:AggHB_complexity}}

When $\mu > 0$, Theorem~\ref{thm:AggHB_main_result} implies \eqref{eq:AggHB_main_result_str_cvx}
    \begin{equation}
        f(\overline{x}_K) - f(x_*) \le \left(1 - \frac{\mu F}{2}\right)^K\frac{4\|x_0 - x_*\|_2^2}{F} \le \frac{4\|x_0-x_*\|^2_2}{F}\mathrm{exp}\left(-\frac{\mu F}{2}K\right).
        \label{eq:AggHB_main_result_str_cvx_2} 
    \end{equation}
Therefore, to ensure that the right-hand side is smaller than $\varepsilon$, number of iterations $K$ should satisfy
    \begin{equation}
        K = \cO\left(\frac{1}{\mu F} \ln \left( \frac{R_0^2}{\varepsilon F} \right)\right),\label{eq:AggHB_compl_tech_cor_1}
    \end{equation}
where $R_0=\|x_0 - x_*\|^2_2$. Assuming $\gamma_i = \gamma$ for $i=1,\ldots,m$ and using constants $\widetilde{\beta}$ and $\hat\beta$ defined as
    \begin{equation*}
        \frac{1}{m}\sum\limits_{i=1}^m \frac{\beta_i}{(1-\beta_i)^2} = \frac{\widetilde{\beta}}{(1-\widetilde{\beta})^2}, \quad \frac{1}{m}\sum\limits_{i=1}^m \frac{1}{1-\beta_i} = \frac{1}{1-\hat{\beta}},
    \end{equation*}
    we get that
    \begin{equation*}
        B = \frac{1}{m}\sum\limits_{i=1}^m \frac{\beta_i \gamma_i \left(1-\beta_i^{K+1}\right)}{(1-\beta_i)^2} \leq \frac{\gamma\widetilde{\beta}}{(1-\widetilde{\beta})^2},\quad F = \frac{1}{m}\sum\limits_{i=1}^{m}\frac{\gamma_i}{1-\beta_i} = \frac{\gamma}{1-\hat\beta}.
    \end{equation*}
    Therefore, the conditions from \eqref{eq:constants_AggHB_2} hold when
    \begin{equation*}
        \frac{\gamma}{1-\hat\beta} \leq \frac{1}{4L},\quad \frac{\gamma^2\widetilde{\beta}}{(1-\widetilde{\beta})^2(1-\hat\beta)} \leq \frac{1 - \max\limits_{i=1,m}\beta_i}{48 L^2}.
    \end{equation*}
    These inequalities together with \eqref{eq:params_AggHB_2} are satisfied for
    \begin{equation}
        \gamma = \min\left\{\frac{\left(1-\max\limits_{i=1,m}\beta_i\right)^2}{2\mu}, \frac{1-\hat{\beta}}{4L}, \frac{(1-\widetilde{\beta})\sqrt{(1-\hat{\beta})\left(1-\max\limits_{i=1,m}\beta_i\right)}}{4\sqrt{3}L\sqrt{\widetilde{\beta}}} \right\}.
        \label{eq:AggHB_compl_tech_cor}
    \end{equation}
Plugging \eqref{eq:AggHB_compl_tech_cor} in \eqref{eq:AggHB_compl_tech_cor_1}, we derive the following complexity result:
    \begin{eqnarray*}
        K &=& \cO \left( \frac{1}{\mu F}\ln \frac{R_0^2}{F \varepsilon}\right) = \cO \left( \frac{1-\hat{\beta}}{\gamma \mu} \ln \frac{(1-\hat{\beta})R_0^2}{\gamma \varepsilon}\right) \\
        &\overset{\eqref{eq:AggHB_compl_tech_cor}}{=}& \cO\Bigg(\max\Bigg(\frac{L}{\mu}, \frac{1-\hat{\beta}}{\left(1-\max\limits_{i=1,m}\beta_i\right)^2}, \frac{L\sqrt{\widetilde{\beta}(1-\hat{\beta})}}{\mu(1-\widetilde{\beta})\sqrt{1-\max\limits_{i=1,m}\beta_i}}\Bigg) \notag\\
        &&\qquad \cdot \ln \Bigg(\frac{R_0^2}{\varepsilon} \cdot\max\Bigg( L, \frac{1-\hat{\beta}}{(1-\max\limits_{i=1,m} \beta_i)^2}, \frac{L \sqrt{\widetilde{\beta}(1-\hat{\beta})}}{(1-\widetilde{\beta})\sqrt{1-\max\limits_{i=1,m} \beta_i}}\Bigg)\Bigg)\\
        &=& \cO\Bigg(\Bigg(\frac{L}{\mu} + \frac{1-\hat{\beta}}{\left(1-\max\limits_{i=1,m}\beta_i\right)^2} + \frac{L\sqrt{\widetilde{\beta}(1-\hat{\beta})}}{\mu(1-\widetilde{\beta})\sqrt{1-\max\limits_{i=1,m}\beta_i}}\Bigg) \notag\\
        &&\qquad \cdot \ln \Bigg(\frac{R_0^2}{\varepsilon} \cdot\Bigg( L + \frac{1-\hat{\beta}}{(1-\max\limits_{i=1,m} \beta_i)^2} + \frac{L \sqrt{\widetilde{\beta}(1-\hat{\beta})}}{(1-\widetilde{\beta})\sqrt{1-\max\limits_{i=1,m} \beta_i}}\Bigg)\Bigg).
    \end{eqnarray*}
    
When $\mu = 0$, Theorem~\ref{thm:AggHB_main_result} implies \eqref{eq:AggHB_main_result_str_cvx}
    \begin{equation}
        f(\overline{x}_K) - f(x_*) \le \frac{4\|x_0 - x_*\|_2^2}{F K}.
        \label{eq:AggHB_main_result_cvx_2} 
    \end{equation}
Therefore, to ensure that the right-hand side is smaller than $\varepsilon$, number of iterations $K$ should satisfy
    \begin{equation}
        K = \cO\left(\frac{R_0^2}{\varepsilon F} \right).\label{eq:AggHB_compl_tech_cor_2}
    \end{equation}
Plugging \eqref{eq:AggHB_compl_tech_cor} in \eqref{eq:AggHB_compl_tech_cor_2}, we get \eqref{eq:AggHB_compl_cvx}.

\end{document}